\documentclass[twoside]{amsart}
\usepackage[centertags]{amsmath}
\usepackage{amsfonts}
\usepackage{dsfont}
\usepackage{amssymb}
\usepackage{amsthm}
\usepackage[ansinew]{inputenc}
\usepackage[cmtip,all]{xy}
\usepackage{graphicx}
\newtheorem{thm}[equation]{Theorem}
\newtheorem{cor}[equation]{Corollary}
\newtheorem{lem}[equation]{Lemma}
\newtheorem{prop}[equation]{Proposition}
\theoremstyle{definition}
\newtheorem{defn}[equation]{Definition}
\theoremstyle{remark}
\newtheorem{rem}[equation]{Remark}

\numberwithin{equation}{section}

\newcommand{\set}[1]{\left\{#1\right\}}

\newcommand{\To}{\longrightarrow}

\newcommand{\C}[1]{\mathbf{#1}} 

\def\op{{op}}
\def\ni{{nil}}
\def\abb{{ab}}

\def\r{\rightarrow} 
\def\rr{\Rightarrow} 

\def\hom{\operatorname{Hom}}
\def\ho{\operatorname{Ho}}
\def\iso{\operatorname{Iso}}
\def\fib{\operatorname{Fib}}

\def\aut{\operatorname{Aut}}

\def\st{\stackrel} 

\newcommand{\ad}{\mathsf{Ad}}

\newcommand{\hopf}{\mathit{Hopf}}

\def\coker{\operatorname{Coker}}
\renewcommand{\ker}{\operatorname{Ker}}

\def\Z{\mathbb{Z}}

\def\N{\mathds{N}}

\def\S{\Sigma}

\def\L{\Omega}

\newcommand{\grupo}[1]{\langle #1\rangle}

\newcommand{\vc}{\Box}    
\newcommand{\vi}{\boxminus}  

\begin{document}

\title[Secondary homotopy groups]{Secondary homotopy groups}%
\author{Hans-Joachim Baues and Fernando Muro}%
\address{Max-Planck-Institut f\"ur Mathematik, Vivatsgasse 7, 53111 Bonn, Germany}%
\email{baues@mpim-bonn.mpg.de, muro@mpim-bonn.mpg.de}%

\thanks{The second author was partially supported
by the project MTM2004-01865 and the MEC postdoctoral fellowship EX2004-0616.}%

\begin{abstract}
Secondary homotopy groups supplement the structure of classical
homotopy groups. They yield a $2$-functor on the groupoid-enriched category of pointed spaces compatible with fiber
sequences, suspensions and loop spaces. They also yield algebraic models of $(n-1)$-connected $(n+1)$-types for
$n\geq 0$.
\end{abstract}
\maketitle

\section*{Introduction}

The computation of homotopy groups of spheres in low degrees in
\cite{toda} uses heavily secondary operations termed Toda brackets.
Such bracket operations are defined by pasting tracks where a track
is a homotopy class of homotopies. Since Toda brackets play a
crucial role in homotopy theory it seems feasible to investigate the
algebraic nature of tracks. Therefore we shift focus from homotopy
groups $\pi_nX$ to secondary homotopy groups
$$\pi_{n,*}X=(\pi_{n,1}X\st{\partial}\To\pi_{n,0}X)$$
defined in this paper. Here $\partial$ is a homomorphism of groups
with $\coker\partial=\pi_nX$ and $\ker\partial=\pi_{n+1}X$, $n\geq1$.

The adjective ``secondary" in the title complements the word ``group", so secondary homotopy groups are \emph{secondary
groups} appearing in homotopy theory. The words ``secondary groups" stand for a variety of dimension $2$
generalizations of the notion of group, like Whitehead's crossed modules \cite{chII}, or reduced or stable
quadratic
modules in the sense of \cite{ch4c}. There is not a well-established terminology in the literature to designate
all these $2$-dimensional ``groups" and we believe that ``secondary groups" has the advantage of being new (so no
confusion with older concepts is created) and short.

The groups $\pi_{n,0}X$ and $\pi_{n,1}X$ are defined directly by use of continuous maps $f\colon S^n\r X$ and
tracks of such maps to the trivial map, so that $\pi_{n,*}X$ is actually a functor in $X$. For $n\geq 2$ the
definition involves the new concept of Hopf invariant for tracks.

We show that the homomorphism $\partial$ has additional algebraic structure, namely $\pi_{1,*}X$
is a crossed module, $\pi_{2,*}X$ is a reduced quadratic module and $\pi_{n,*}X$, $n\geq 3$, is a stable
quadratic module.

Crossed modules were introduced by J. H. C. Whitehead in \cite{chII} and, in fact, for a reduced $CW$-complex $X$
with $1$-skeleton $X^1$ our
secondary homotopy group $\pi_{1,*}X$ is weakly equivalent to the crossed module
$$\pi_2(X,X^1)\To\pi_1X^1$$ studied by \cite{chII}. Similarly if $X$ is an $(n-1)$-reduced $CW$-complex then $\pi_{n,*}X$ for $n\geq 2$ is weakly equivalent to
the quadratic modules obtained in \cite{ch4c} in terms of the cell structure of $X$. 

The topological and functorial definition of secondary homotopy groups $\pi_{n,*}X$ in this paper is crucial to understand
new properties of the corresponding concepts in the literature. For example, we are able to determine the algebraic
properties of the loop and suspension operators on secondary homotopy groups. As main new results, we describe
the fiber sequence for secondary homotopy groups, and we show that secondary homotopy groups form a $2$-functor
on the groupoid-enriched category of pointed spaces. We also prove that secondary homotopy groups are algebraic models for
two-stage spaces, i.e. spaces with only two non-trivial homotopy groups in consecutive dimensions. Such two-stage
spaces generalize Eilenberg-MacLane spaces, which are spaces with a single non-trivial homotopy group.


The crucial topological tool that we introduce in this paper for the construction of the secondary homotopy
groups is the Hopf invariant for tracks. This very basic piece of homotopy theory does not seem to have been
considered before. It generalizes the classical Hopf invariant for maps. In this paper we only need to compute
Hopf invariants for tracks between one-point unions of spheres. However we believe that a notion of Hopf invariant for more general
tracks would deserve to be studied elsewhere.

Hopf invariant computations for tracks are the basis of the applications obtained in the sequels of this paper,
see \cite{2hg2,2hg3,2hg4}.
These computations are achieved
by using geometric tools such as orthogonal group actions and Clifford algebras. This stresses the relevance of
the definition of secondary homotopy groups by topological means.

In \cite{2hg3} we determine the algebraic nature of smash
product operations on secondary homotopy groups. For the (stable)
secondary homotopy groups of spectra this leads in \cite{2hg4} to an algebraic
``tensor product'' approximating the smash product of spectra on a secondary level.

The computation of the algebra of secondary cohomology operations in
\cite{asco} shows examples where secondary homotopy groups can be
algebraically determined successfully. It is the aim of the authors to generalize the theory of \cite{asco},
concerning the Eilenberg-MacLane spectrum, for
general spectra.

Moreover, we will discuss in a sequel of this paper generalized Whitehead products for secondary homotopy
groups. In fact, J. H. C. Whitehead introduced in \cite{arhg}
Whitehead products as an additional algebraic structure of homotopy
groups. We may consider the secondary homotopy groups together with their algebraic properties also as such
an enriching structure.

\subsection*{Further connections with the literature}

The secondary groups considered in this paper (i.e. crossed modules and reduced or stable quadratic modules) are
equivalent, at least up to homotopy, to some other algebraic structures in the literature, such as categorical
groups, which can also be braided or symmetric. There are strict and non-strict categorical groups. Strict algebraic
objects, like the secondary groups we use, are often more convenient to work with. Garz\'on, Miranda and del R\'io
\cite{tehgte} endow the fundamenthal groupoid of an iterated loop space with the structure of a categorical
group, which can be braided or symmetric according to the number of loopings. Our secondary homotopy groups could
be regarded as strictifications of these. Indeed the objects of the categorical groups in \cite{tehgte} are the
same as the generators of our secondary homotopy groups. The connection between  the secondary groups that we use and categorical
groups can be established through Conduch\'e's $2$-crossed modules \cite{2cm}, see \cite{ccscg}, \cite{3dnac} and
\cite{2tils}.

More precisely, the equivalence between categorical groups and $2$-crossed modules can be found in \cite{ccscg},
\cite{3dnac} where the authors define strict ``homotopy categorical groups'' for simplicial groups. Given a pointed
space $X$ one can consider the Kan loop group $GS_\bullet X$ on the reduced singular simplicial set of the pointed
space $X$. The categorical groups defined in \cite{ccscg},
\cite{3dnac} applied to $GS_\bullet X$ are models for two-stage spaces. Since the same result holds for our
secondary homotopy groups there must be a ``weak equivalence'' between both constructions. The definition of an
explicit equivalence is however out of the scope of this paper.

The connection beween $2$-crossed modules and secondary groups is established in \cite{2tils}, where the authors
show that quadratic modules are highly connected $2$-crossed modules. Indeed they have nilpotency degree $2$, so
they are just one step away from being abelian. This yields an advantage of secondary homotopy groups over the
categorical groups constructed in \cite{ccscg} and \cite{3dnac}. In addition our construction does not rely on
the Kan loop group $GS_\bullet X$. On the contrary it is purely topological since it uses continuous maps and
tracks. 

\subsection*{Acknowledgement}

The authors are very grateful to the referee for his careful reading of the paper and for many valuable
suggestions.




\section{Tracks between maps}\label{tbm}


We consider the category $\C{Top}^*$ of compactly generated pointed spaces $X=(X,*)$ and pointed maps $f\colon X\r Y$. For any (unpointed) space $X$ we define $X_+=X\sqcup\set{*}$ as
the same space with an outer base-point $*$. The smash product of two pointed spaces is defined by
$$X\wedge Y=(X\times Y)/(X\times *\cup *\times Y).$$
It is associative and commutative, since it defines a symmetric monoidal structure on the category $\C{Top}^*$.

Homotopies $IX\r Y$ are defined by using the
reduced cylinder $IX=I_+\wedge X$, where $I=[0,1]$ is the unit interval, with structure maps
\begin{equation}
X\vee X\st{i}\To IX\st{p}\To X.
\end{equation} Here $\vee$ is the symbol for the coproduct, $i$ is the inclusion of the boundary and $p$ is the projection.
Given two maps $f,g\colon X\r Y$ a \emph{track} $H\colon f\rr g$ is a homotopy class of homotopies $IX\r Y$, 
from $f$ to $g$, relative to the boundary.
By abuse of language we denote a homotopy and the represented track by the same symbol.
 In diagrams tracks will be denoted as follows.
\begin{equation}\label{track}
\xymatrix@C=35pt{X\ar@/^10pt/[r]^f_{\;}="a"\ar@/_10pt/[r]_{g}^{\;}="b"&
Y.\ar@{=>}"a";"b"^{H}}
\end{equation}
The trivial track $0^\vc_f\colon f\rr f$ is represented by $fp\colon IX\r Y$ and the inverse of a track 
$H\colon f\rr g$ is $H^\vi\colon g\rr f$.
The \emph{vertical composition} of
tracks
$$\xymatrix@C=40pt{X\ar@/^20pt/[r]^f_{\;}="a"\ar[r]|{g}^{\;}="b"_{\;}="c"\ar@/_20pt/[r]_h^{\;}="d"&Y\ar@{=>}"a";"b"^H\ar@{=>}"c";"d"^K}$$
is defined by pasting homotopies representing $H$ and $K$ and is denoted by 
\begin{equation*}
\xymatrix@C=50pt{X\ar@/^15pt/[r]^f_{\;}="a"\ar@/_15pt/[r]_{h}^{\;}="b"&
Y\ar@{=>}"a";"b"^{K\vc H}}
\end{equation*}
One can also compose horizontally a track as in diagram
(\ref{track}) with maps $k\colon W\r X$ and $l\colon Y\r Z$ to
obtain tracks
$$Hk\colon fk\rr gk\text{   and   }lH\colon lf\rr lg$$
in the obvious way.
If we have a diagram like
\begin{equation*}
\xymatrix@C=35pt{X\ar@/_10pt/[r]_g^{\;}="a"\ar@/^10pt/[r]^{f}_{\;}="b"&Y
\ar@/_10pt/[r]_{g'}^{\;}="c"\ar@/^10pt/[r]^{f'}_{\;}="d"&Z\ar@{=>}"b";"a"^H\ar@{=>}"d";"c"^{H'}}
\end{equation*}
the equality
$$(g'H)\vc(H'f)=(H'g)\vc(f'H)$$
holds and this element is the \emph{horizontal composition} of $H$ and $H'$ denoted by juxtaposition
$$\xymatrix@C=50pt{X\ar@/^15pt/[r]^{f'f}_{\;}="a"\ar@/_15pt/[r]_{g'g}^{\;}="b"&
Z\ar@{=>}"a";"b"^{H'H}}$$
Tracks endow $\C{Top}^*$ with the structure of a groupoid-enriched category. 

Recall that a \emph{groupoid-enriched category} $\C{C}$ is the same as a $2$-category where all $2$-morphisms are
vertically invertible. In this paper we will work with some other groupoid-enriched categories, apart from
$\C{Top}^*$. We will always use the same notation as above for $2$-morphisms, 
identity $2$-morphisms (also called trivial $2$-morphisms), vertical and horizontal composition, and vertical
inverses. A morphism $f\colon X\r Y$ in $\C{C}$ has an automorphism group $\aut_\C{C}(f)$ given by the set of
$2$-morphisms
$f\rr f$  and the vertical composition. The \emph{homotopy category} of $\C{C}$ is the ordinary category $\C{C}/\!\simeq$ obtained
by identifying two morphisms $f,g\colon X\r Y$ in $\C{C}$ provided there exists a $2$-morphism between them
$H\colon f\rr g$. For instance if $\C{C}=\C{Top}^*$ then the homotopy category is the usual one. A $2$-functor $\varphi\colon\C{C}\r\C{D}$ between groupoid-enriched categories with the same
objects which is the identity on objects is said to be a \emph{weak equivalence} if $\varphi$ induces an isomorphism
between the automorphism groups $\aut_\C{C}(f)\cong\aut_\C{D}(\varphi(f))$ for any morphism $f$ in $\C{C}$ and an
isomorphism between the homotopy categories $\C{C}/\!\simeq\;\cong\;\C{D}/\!\simeq$.

\begin{rem}\label{*}
The groupoid-enriched category $\C{Top}^*$ has a strict zero object, the one-point space $*$. In particular zero maps are
defined. Such a groupoid-enriched category has the crucial property that any $2$-morphism composed with a zero map becomes
automatically an identity $2$-morphism.
\end{rem}

Maps from a coproduct $X\vee Y$ in $\C{Top}^*$ are given by pairs of maps $(f_1,f_2)\colon X\vee Y\r Z$. Similarly a track $H\colon (f_1,f_2)\rr(g_1,g_2)$
between maps $(f_1,f_2),(g_1,g_2)\colon X\vee Y\r Z$ is given by a pair of tracks $H=(H_1,H_2)$ with $H_i\colon
f_i\rr g_i$ $(i=1,2)$. This means that the one-point union of pointed spaces is an enriched coproduct in
$\C{Top}^*$.

The suspension $\Sigma X$ is the quotient space $IX/(X\vee X)=S^1\wedge X$. It defines a $2$-functor
$$\S\colon\C{Top}^*\To\C{Top}^*.$$
The suspension of a track $H\colon f\rr g$ between maps $f,g\colon X\r Y$ represented by a homotopy $H\colon IX\r
Y$ is the track $\S H\colon \S
f\rr\S g$ represented by the homotopy
$$I\S X=I_+\wedge S^1\wedge X\cong S^1\wedge I_+\wedge X\st{S^1\wedge H}\To S^1\wedge Y=\S Y.$$
We will use the identifications
\begin{equation}\label{idsus}
\Sigma (X\vee Y)=\Sigma X\vee\Sigma Y,\;\;\;\Sigma^nS^0= S^{n}, \;\; n\geq 0.
\end{equation}
For the definition of homotopy groups we choose a particular co-H-group structure on $S^1$ given by maps $\mu\colon S^1\r S^1\vee S^1$ and $\nu\colon S^1\r S^1$
satisfying the usual properties. We use explicitly these maps in many constructions throughout this paper,
however these constructions do not depend
on this choice since the maps $\mu$ and $\nu$ are unique up to a canonical track.

The loop space functor $\L$ is the right-adjoint of the suspension $\S$. The adjoint of a map $f\colon\S X\r Y$ is denoted by $ad(f)\colon X\r\L Y$. The adjoint of the identity map $1\colon\S X\r \S X$
is a natural inclusion
\begin{equation}\label{ad1}
ad(1)\colon X\hookrightarrow \L\S X.
\end{equation}
As a pointed set the $n$-fold loop space $\L^n X$ is the set of
pointed maps $S^n\r X$ and the base-point corresponds to the trivial map.
By using the interchange homeomorphism of the smash product we see that suspensions and cylinders commute up to
natural isomorphism in $\C{Top}^*$, $I\S X\cong\S I X$.
However one has to be careful with signs because the interchange of factors in $S^1\wedge S^1=S^2$ induces $-1$ on the homotopy group $\pi_2$.

\section{Groups of nilpotency degree $2$}

Consider the forgetful functor from groups to pointed sets
$\C{Gr}\To\C{Set}^*.$
This functor has a left adjoint
$$\grupo{\cdot}\colon\C{Set}^*\To\C{Gr}\colon A\mapsto \grupo{A}.$$
Here $\grupo{A}$ is the quotient of the free group with basis $A$ by the normal subgroup generated by the base-point $*\in A$. This group
is isomorphic to the free group with basis $A-\set{*}$.
We denote $\vee_AS^1=\S A$.
As usual we identify the fundamental group
of $\vee_AS^1$ with a free group, i.e. $\pi_1(\vee_AS^1)=\grupo{A}$.
The free group of nilpotency class $2$ (\emph{free nil-group} for short), generated by the pointed set $A$, is the quotient
$$\grupo{A}_\ni=\frac{\grupo{A}}{\Gamma_3\grupo{A}}$$
where $\Gamma_3\grupo{A}$ is the $3^\text{rd}$ term of the lower central series of $\grupo{A}$, i.e. the
subgroup generated by triple commutators 
$[x,[y,z]]$ $(x,y,z\in\grupo{A})$.  In this paper we always write group laws additively, even for non-abelian
groups, so that the commutator is $[x,y]=-x-y+x+y$.
The free abelian group $\Z[A]$ on a pointed set $A$ is the abelianization of $\grupo{A}$ and of $\grupo{A}_\ni$.

If $\C{gr}$, $\C{nil}$ and $\C{ab}$ are the categories of free groups, free nil-groups and free abelian groups, respectively, then there are
obvious nilization and abelianization functors
\begin{equation}\label{abnil}
\begin{array}{cc}
\xymatrix{\C{gr}\ar[rr]^\abb\ar[rd]_\ni&&\C{ab}\\&\C{nil}\ar[ru]_\abb&}&
\xymatrix{\grupo{A}\ar@{|->}[rr]\ar@{|->}[rd]&&\Z[A]\\&\grupo{A}_\ni\ar@{|->}[ru]&}
\end{array}
\end{equation}

Let $\grupo{A}_\ni\twoheadrightarrow\Z[A]$ be the natural projection
carrying $x$ to $\set{x}$. Since the commutator bracket in $\grupo{A}_\ni$
is bilinear the homomorphism
\begin{equation}\label{kio}
\partial\colon\otimes^2\Z[A]\r\grupo{A}_\ni,\;\;\;
\partial(\set{x}\otimes \set{y})=[x,y],
\end{equation}
is well defined. Here the tensor
square of an abelian group $A$ is denoted by $\otimes^2A=A\otimes
A$. This fact is relevant since it allows us to define a groupoid enrichment of $\C{nil}$ which is connected to
the tracks between one-point unions of spheres through the Hopf invariant, see the next section. This connection and its consequences
are crucial for the definition of secondary homotopy groups.

Let $T\colon A\otimes B\r B\otimes A$ be the interchange
isomorphism $T(a\otimes b)=b\otimes a$. 
The reduced tensor square is
the following cokernel
$$\otimes^2A\st{1+T}\To\otimes^2A\st{\bar{\sigma}}\twoheadrightarrow\hat{\otimes}^2A.$$
We denote $\bar{\sigma}(a\otimes b)=a\hat{\otimes}b$.
We define the functor $\otimes^2_n$ as
$$\otimes^2_n=\left\{\begin{array}{ll}\otimes^2,&\text{if $n=2$;}\\&\\\hat{\otimes}^2,&
\text{if $n\geq 3$.}\end{array}\right.$$
Here we write  $a\otimes b\in\otimes^2_nA$ with $a\otimes b=a\hat{\otimes}b$ for $n\geq 3$.
Moreover, $\Gamma_n$ is the functor
$$\Gamma_n=\left\{\begin{array}{ll}\Gamma,&\text{if $n=2$;}\\&\\-\otimes\Z/2,&
\text{if $n\geq 3$;}\end{array}\right.$$
where $-\otimes\Z/2$ is the ordinary tensor product of abelian groups and $\Gamma$ is Whitehead's universal quadratic functor, see \cite{aces}.
There is a natural exact sequence
\begin{equation}\label{exagam}
\Gamma_n\Z[A]\hookrightarrow\otimes^2_n\Z[A]\st{\partial}\To\grupo{A}_\ni\twoheadrightarrow\Z[A].
\end{equation}
Here the first arrow is induced by the function sending $x\in\Z[A]$
to $x\otimes x\in\otimes^2_n\Z[A]$, see for example \cite{ch4c}. Moreover,
these exact sequences fit into a natural commutative diagram
$$\xymatrix{\Gamma\Z[A]\ar@{->>}[d]_\sigma\ar@{^{(}->}[r]&\otimes^2\Z[A]\ar[r]^\partial\ar@{->>}[d]_{\bar{\sigma}}&\grupo{A}_\ni\ar@{=}[d]\ar@{->>}[r]&\Z[A]\ar@{=}[d]\\
\Z[A]\otimes\Z/2\ar@{^{(}->}[r]&\hat{\otimes}^2\Z[A]\ar[r]^\partial&\grupo{A}_\ni\ar@{->>}[r]&\Z[A]}$$
where the upper $\partial$ is the homomorphism in (\ref{kio}) above, and the lower $\partial$ is a factorization
of (\ref{kio}) through the reduced tensor square.

\section{Nil-tracks and Hopf invariants of tracks}

In this section we introduce a new homotopy invariant which is crucial for the definition of secondary homotopy
groups and for the further development of this theory. It is a Hopf invariant for tracks which generalizes the most
classical Hopf invariant for maps, see below. In this paper we only consider the Hopf invariant of tracks between
one-point unions of 
spheres. However it would be interesting to develop a more general notion.

\begin{defn}\label{niltrack}
Let $f,g$ be maps $f,g\colon S^1\r\vee_AS^1$ where $A$ is a discrete pointed set, and let
$$\S^{n-1}f, \S^{n-1} g\colon S^n\r\vee_AS^n$$ be their $(n-1)$-fold
suspensions, $n\geq 1$. A track $H\colon\S^{n-1}f\rr\S^{n-1}g$, represented by a
homotopy $H\colon IS^n\r \vee_AS^n$, is said to be a \emph{nil-track}
if the adjoint
$$ad(H)\colon IS^1\To\L^{n-1}(\vee_AS^n)$$
of the map
$$S^{n-1}\wedge I_+\wedge S^1\cong I_+\wedge S^{n-1}\wedge S^1=IS^n\st{H}\To \vee_AS^n.$$
induces a trivial homomorphism
$$0=H_2ad(H)\colon H_2(I S^1, S^1\vee S^1)\To H_2(\L^{n-1}(\vee_A S^n),\vee_AS^1).$$
The adjoint of $H$ sends the boundary of the cylinder $IS^1$ into $\vee_AS^1$ since $H$ restricted to the boundary
is an $(n-1)$-fold suspension. Of course for $n=1$ all tracks $H$ above are nil-tracks since $H_2ad(H)$ maps to the trivial group.

Let $f,g$ be now maps between wedges of $1$-spheres $f,g\colon \vee_B S^1\r\vee_AS^1$, and let $\S^{n-1}f, \S^{n-1} g$ be their $(n-1)$-fold suspensions. A track
$H\colon \S^{n-1}f\rr \S^{n-1}g$ is a \emph{nil-track} if all restricted tracks $Hi_b$ are nil-tracks where $i_b\colon
S^n\r\vee_BS^n$ is the inclusion given by $b\in B-\set{*}$.
\end{defn}

The homology groups involved in the definition of nil-tracks are computable. Indeed,
$$H_2(I S^1, S^1\vee S^1)\cong H_2(\S S^1)=H_2S^2\cong\Z.$$
Moreover,
$$H_2(\L^{n-1}(\vee_AS^n))\st{\cong}\To H_2(\L^{n-1}(\vee_A S^n),\vee_AS^1)$$
is an isomorphism, and the Pontrjagin product
$$\otimes^2\Z[A]=H_1(\L^{n-1}(\vee_AS^n))\otimes H_1(\L^{n-1}(\vee_AS^n))\To H_2(\L^{n-1}(\vee_AS^n))$$
is an isomorphism for $n=2$ and induces an isomorphism for $n\geq 2$
\begin{equation}\label{uniso}
{\otimes}^2_n\Z[A]\cong H_2(\L^{n-1}(\vee_AS^n)),
\end{equation}
compare notation in (\ref{exagam}).

\begin{defn}\label{hopf}
Let $n\geq 2$. Given a track $H\colon\S^{n-1}f\rr\S^{n-1}g$ for maps $f,g\colon
S^1\r\vee_AS^1$ the \emph{Hopf invariant} of $H$ is defined as
$$\hopf(H)=(H_2 ad(H))(1)\in\otimes^2_n\Z[A],$$
where we apply the homology functor $H_2$  as in Definition \ref{niltrack}.
In particular, $H$ is a nil-track if and only if $\hopf(H)=0$. More
generally, if $H\colon\S^{n-1}f\rr\S^{n-1}g$ is a track for maps
$f,g\colon \vee_BS^1\r\vee_AS^1$ the Hopf invariant of $H$ is the
homomorphism
$$\hopf(H)\colon\Z[B]\To\otimes^2_n\Z[A]$$
defined by $\hopf(H)(b)=\hopf(Hi_b)$, where $i_b\colon
S^1\subset\vee_B S^1$ is the inclusion of the factor corresponding
to $b\in B-\set{*}$. Such a track $H$ is a nil-track if and only if
$\hopf(H)=0$.

In case $n=1$ we have $\hopf(H)=0$ for any track $H$ as above.
\end{defn}

\begin{rem}\label{conclasi}
Any element $x\in\pi_{3}\vee_AS^2$ determines a track $x\colon 0\rr 0$ for the trivial map $0\colon S^2\r \vee_AS^2$. This track is given by
the homotopy $IS^2\r\S S^2=S^3\st{x}\r \vee_AS^2$, where the first map is the obvious projection.
The reader can check that $-\hopf(x)$ is the classical Hopf invariant of $x$. The sign is due to the fact that in order to define the Hopf
invariant of $x$ as a track we need to consider the map $I_+\wedge S^1\wedge S^1\cong S^1\wedge I_+\wedge S^1$
interchanging the first two factors of the smash product and this map induces
$-1\colon S^3\r S^3$ up to homotopy.
\end{rem}


The next results are crucial for this paper.

\begin{thm}\label{niltrackes}
Let $f,g\colon \vee_AS^1\r\vee_BS^1$ be two maps and $n\geq 1$. If a nil-track
$$N_{f,g}\colon\S^{n-1}f\rr\S^{n-1}g$$
exists then it is unique. Moreover, $N_{f,g}$ exists if and only if
\begin{itemize}
\item $\pi_1f=\pi_1g\colon\grupo{A}\r\grupo{B}$, if $n=1$;
\item or $(\pi_1f)_\ni=(\pi_1g)_\ni\colon\grupo{A}_\ni\r\grupo{B}_\ni$, if $n\geq 2$.
\end{itemize}
Furthermore, trivial tracks are nil-tracks and the vertical and
horizontal composition of nil-tracks are also nil-tracks.
\end{thm}

The case $n=1$ is obvious since a one-point union of circles is aspherical, so a track between maps $f,g\colon
\vee_AS^1\r \vee_BS^1$ is necessarily unique provided it exists, it is always a nil-track as we have noticed
above, and it is well-known that $H$ exists if and only if $f$ and $g$ induce the same homomorphism on
fundamental groups.
For $n\geq 2$ this theorem is a immediate consequence of the following one. 

\begin{thm}\label{propi}
Let $n\geq 2$ and let $f,g\colon\vee_AS^1\r\vee_BS^1$ be maps such that for any $x\in\grupo{A}_\ni$ we have
$(\pi_1g)_\ni(x)=(\pi_1f)_\ni(x)+\partial\alpha(x)$
for some homomorphism $\alpha\colon\Z[A]\r\otimes^2_n\Z[B]$. Then 
there exists a unique track $H\colon\S^{n-1}f\rr\S^{n-1}g$ with Hopf invariant
$\hopf(H)=\alpha$ and conversely. Moreover,
the Hopf invariant of tracks satisfies the following formulas. Given
a diagram
$$\xymatrix@C=40pt{\vee_AS^n\ar@/^30pt/[r]^{\S^{n-1}f}_{\;}="a"
\ar[r]|{\S^{n-1}g}^{\;}="b"_{\;}="c"\ar@/_30pt/[r]_{\S^{n-1}h}^{\;}="d"&\vee_BS^n\ar@{=>}"a";"b"^H\ar@{=>}"c";"d"^K}$$
the equation
\begin{enumerate}
\item $\hopf(K\vc H)=\hopf(K)+\hopf(H)$ holds.
\end{enumerate}
Furthermore, if we consider the diagram
\begin{equation*}
\xymatrix@C=35pt{\vee_AS^n\ar[r]^{\S^{n-1}k}&\vee_BS^n\ar@/_15pt/[r]_{\S^{n-1}g}^{\;}="a"
\ar@/^15pt/[r]^{\S^{n-1}f}_{\;}="b"&\vee_CS^n
\ar[r]^{\S^{n-1}h}&\vee_DS^n\ar@{=>}"b";"a"^H}
\end{equation*}
then
\begin{enumerate}\setcounter{enumi}{1}
\item $\hopf(H(\S^{n-1}k))=\hopf(H)(\pi_1k)_\abb$,
\item $\hopf((\S^{n-1}h)H)=(\otimes^2_n(\pi_1h)_\abb)\hopf(H)$.
\end{enumerate}
In addition given a track $H\colon\S^{n-1}f\rr\S^{n-1}g$ between maps $f,g,\colon\vee_AS^1\r\vee_BS^1$ one gets the following equations.
\begin{enumerate}\setcounter{enumi}{3}
\item $\hopf(\S H)=0$ if $n=1$,
\item $\hopf(\S H)=\bar{\sigma}\hopf(H)$ if $n=2$,
\item $\hopf(\S H)=\hopf(H)$ if $n\geq 3$.
\end{enumerate}
\end{thm}

This theorem is a simple consequence of the theory developed in \cite{ch4c} that we now recall.

Let $\C{S}(n)\subset\C{Top}^*$ be the full groupoid-enriched subcategory of
one-point unions of $n$-spheres and let $\C{gr}$ be the
category of free groups regarded as a groupoid-enriched
category with only the trivial $2$-morphisms. Then there is a $2$-functor
$$\pi_1\colon\C{S}(1)\To\C{gr}$$ given by the fundamental group, $\pi_1(\vee_AS^1)=\grupo{A}$.
This $2$-functor is a weak equivalence. This follows easily from
\cite{ch4c} VI.3.13 and the fact that one-point unions of $1$-spheres do not
have higher-dimensional homotopy groups.

For $n\geq 2$ we consider the groupoid-enriched subcategory
$\bar{\C{S}}(n)\subset\C{S}(n)$ of suspended maps. Here objects of $\bar{\C{S}}(n)$ are
one-point unions of $1$-spheres $\vee_AS^1$, maps $f,g\colon
\vee_AS^1\r\vee_BS^1$ in $\bar{\C{S}}(n)$ are maps in $\C{Top}^*$ and $2$-morphisms $H\colon
f\rr g$ in $\bar{\C{S}}(n)$ are tracks $H\colon\S^{n-1}f\rr\S^{n-1}g$ in $\C{Top}^*$.
The inclusion
\begin{equation}\label{fuli}
\S^{n-1}\colon\bar{\C{S}}(n)\subset\C{S}(n)
\end{equation}
is given by the $(n-1)$-fold suspension on objects and morphisms and
it is the identity on $2$-morphisms. This is actually a weak equivalence of
groupoid-enriched categories. See \cite{ch4c} VI.4.7.

We now consider the algebraic goupoid-enriched category $\C{nil}(n)$ defined as
follows. Objects and morphisms are the same as in $\C{nil}$. A $2$-morphism
$\alpha\colon \varphi\rr\psi$ between homomorphisms
$\varphi,\psi\colon\grupo{A}_\ni\r\grupo{B}_\ni$ is a homomorphism
$\alpha\colon \Z[A]\r\otimes^2_n\Z[B]$ such that
$\varphi(x)+\partial\alpha (\set{x})=\psi(x)$ for any
$x\in\grupo{A}_\ni$, i.e. a $2$-morphism $\alpha\colon\varphi\rr\psi$ can only exist if the abelianizations coincide
$\varphi^\abb=\psi^\abb$, and in this case $\alpha$ is a lift of the pointwise difference $-\varphi(x)+\psi(x)$,
which lies in the commutator subgroup of $\grupo{B}_\ni$, to the tensor square $\otimes^2\Z[B]$, compare the
exact sequence (\ref{exagam}). The vertical composition is given by addition
of abelian group homomorphisms, and for the horizontal composition
one uses the abelianization functor $ab\colon\C{nil}\r\C{ab}$ and
the bifunctor
$$\hom_\C{ab}(-, \otimes^2_n)\colon\C{ab}^\op\times\C{ab}\To\C{Ab}.$$


For any $n\geq 2$ there is a weak equivalence of groupoid-enriched categories
$$\hopf\colon\bar{\C{S}}(n)\To\C{nil}(n)$$ defined by
$\vee_AS^1\mapsto\grupo{A}_\ni$, $f\mapsto(\pi_1f)_\ni$ and
$H\mapsto\hopf(H)$, where we use the Hopf invariant for tracks. Compare \cite{ch4c} VI.4.7.
This weak equivalence is compatible on the left hand side with the suspension functor
$$\S\colon\bar{\C{S}}(n)\To\bar{\C{S}}(n+1)$$
which is the identity on objects and morphisms and on tracks it is given by the suspension of tracks in 
$\C{Top}^*$, and on the right hand side with the $2$-functors
$$\C{gr}\To\C{nil}(2)\To\C{nil}(n),\;\;n\geq3,$$
given by the nilization and the natural projection
$\bar{\sigma}\colon\otimes^2\twoheadrightarrow\hat{\otimes}^2$ respectively. Here we set $\bar{\C{S}}(1)=\C{S}(1)$ for $n=1$. 
Theorem
\ref{propi}, and therefore Theorem \ref{niltrackes}, follows readily from this. 

\section{Secondary homotopy groups of a pointed space}

We now introduce secondary homotopy groups which enrich the
structure of the classical homotopy groups $\pi_nX$ of a pointed space.

\begin{defn}\label{2hg}
Let $X$ be a pointed space and $n\geq 1$. The \emph{secondary homotopy group} $\pi_{n,*}X$ is the map
$$\partial\colon\pi_{n,1}X\To\pi_{n,0}X$$
defined as follows. 
Let
$$\pi_{n,0}X=\left\{\begin{array}{ll}\grupo{\L X},&n=1;\\{}\\\grupo{\L^n X}_\ni,&n\geq 2.\end{array}\right.$$
Here the $n$-fold loop space, regarded as a discrete pointed set, generates a free (nil-)\\group. Moreover, $\pi_{n,1}X$ is the set of equivalence classes $[f,F]$ 
represented
by a map $f\colon S^1\r \vee_{\L^n X}S^1$ and a track
$$\xymatrix{S^n\ar[r]_{\S^{n-1}f}^<(.98){\;\;\;\;\;}="a"\ar@/^25pt/[rr]^0_{}="b"&S^n_X\ar[r]_{ev}&X\ar@{=>}"a";"b"_F}$$
Here the pointed space $$S^n_X=\vee_{\L^n X}S^n=\S^n\L^n X$$ is the $n$-fold suspension of the $n$-fold loop space
$\L^nX$, where $\L^nX$ is regarded as a pointed set with the discrete topology. Hence $S^n_X$ is
the coproduct of $n$-spheres indexed by the set of non-trivial maps $S^n\r X$, and $ev\colon S^n_X\r X$
is the obvious evaluation map. Moreover, for the sake of simplicity given a map
$f\colon S^1\r \vee_{\L^n X}S^1$ we will denote $f_{ev}=ev(\S^{n-1}f)$, so that $F$ in the previous diagram is a track $F\colon f_{ev}\rr 0$.
The equivalence relation $[f,F]=[g,G]$ holds provided the nil-track $N_{f,g}\colon\Sigma^{n-1}f\rr\S^{n-1}g$ exists, see Theorem \ref{niltrackes}, and
the composite track in the following diagram is the trivial track.
$$\xymatrix@C=50pt{S^n\ar@/^40pt/[rr]^0_{\;}="a"\ar@/_40pt/[rr]_0^{\;}="f"\ar@/^15pt/[r]|{\S^{n-1}f}^<(.935){\;}="b"_{\;}="c"\ar@/_15pt/[r]|{\S^{n-1}g}^{\;}="d"_<(.93){\;}="e"
&S^n_X\ar[r]^{ev}&X\ar@{=>}"a";"b"^{F^\vi}\ar@{=>}"c";"d"^{N_{f,g}}\ar@{=>}"e";"f"^G}$$
That is $F=G\vc(ev\,N_{f,g})$. The map
$\partial$ is defined by the  formula
$$\partial[f,F]=\left\{\begin{array}{ll}(\pi_1f)(1),&n=1,\\{}\\(\pi_1 f)_\ni(1)       ,&n\geq 2,\end{array}\right.$$
where $1\in\pi_1S^1=\Z$.

A map $g\colon X\r Y$ in $\C{Top}^*$ induces a map
$\pi_{n,*}g\colon\pi_{n,*}X\r\pi_{n,*}Y$ given by the following
commutative diagram.
\begin{equation}\label{inmor}
\xymatrix{\pi_{n,1}X\ar[r]^{\pi_{n,1}g}\ar[d]_\partial&\pi_{n,1}Y\ar[d]^\partial\\\pi_{n,0}X\ar[r]_{\pi_{n,0}g}&\pi_{n,0}Y}
\end{equation}
Here the lower homomorphism $\pi_{n,0}g=\grupo{\L^n g}$ is induced
by the map of pointed sets $\L^ng\colon\L^nX\r\L^nY$. Moreover, an
element $[f,F]\in\pi_{n,1}X$ is sent by the upper map $\pi_{n,1}g$
to $[(\S\L^n g) f,gF]$,
$$\xymatrix{&&X\ar@/^10pt/[rd]^g&\\S^n\ar@/^15pt/[rru]^0_{\;}="b"\ar[r]_{\S^{n-1}f}&S^n_X\ar[ru]_{ev}^{\;}="a"
\ar[r]_{\S^n\L^ng}&S^n_Y\ar[r]_{ev}&Y\ar@{=>}"a";"b"^F}$$
\end{defn}

We also define $\pi_{n,*}X$ for $n=0$ as follows.

\begin{defn}
For $n=0$ let $\pi_{0,*}X$ be the \emph{fundamental pointed groupoid} of the pointed space $X$ for which $\pi_{0,0}X$ 
is  $X$ regarded as a discrete pointed set 
and $\pi_{0,1}X$ is the set of tracks between points in $X$. For this we
recall that a \emph{pointed groupoid} is a small category $\C{G}$
with a distinguished object $*\in Ob\C{G}$ such that all morphisms
are isomorphisms. A morphism of pointed groupoids
$F\colon\C{G}\r\C{H}$ is a functor preserving the distinguished
object $F(*)=*$, and the category of pointed groupoids is denoted by
$\C{grd}^*$. The morphism $F$ is a \emph{weak equivalence} if it
induces a bijection between the pointed sets of isomorphism classes
of objects $\iso(F)\colon\iso(\C{G})\cong\iso(\C{H})$ and if
$F\colon\aut_\C{G}(x)\cong\aut_\C{H}(F(x))$ is an isomorphism for
any object $x$ in $\C{G}$. The fundamental groupoid is a functor
$\pi_{0,*}\colon\C{Top}^*\r\C{grd}^*$ in the obvious way. 
\end{defn}

We now study the algebraic structure of secondary homotopy groups $\pi_{n,*}X$ with $n\geq 1$.

\begin{prop}\label{grupoes}
For $n\geq 1$ there is a group structure on $\pi_{n,1}X$ such that the map
$\partial\colon\pi_{n,1}X\r\pi_{n,0}X$ is a homomorphism. Moreover, the rows in (\ref{inmor}) are also group
homomorphisms.
\end{prop}

\begin{proof}
We define
the sum of two elements $[f,F],[g,G]\in\pi_{n,1}X$ by the following diagram
$$\xymatrix{S^n\ar[r]^(.4){\Sigma^{n-1}\mu}&S^n\vee S^n\ar[rr]_-{(\S^{n-1}f,\S^{n-1}g)}^<(.66){\;\;\;\;\;}="a"\ar@/^30pt/[rrr]^0_{\;}="b"&&
S^n_X\ar[r]_{ev}&X\ar@{=>}"a";"b"_{(F,G)}}$$
i.e.
$$[f,F]+[g,G]=[(f,g)\mu,(F,G)(\S^{n-1}\mu)].$$
One can readily check by using Theorem \ref{niltrackes} that this
operation 
is associative
and $[0,0^\vc]$ is a unit element. 
The inverse of an element
$[f,F]$ is represented by
$$\xymatrix{S^n\ar[r]_{\S^{n-1}\nu}&S^n\ar[r]_{\S^{n-1}f}^<(.98){\;\;\;\;\;}="a"\ar@/^25pt/[rr]^0_{}="b"&S^n_X\ar[r]_{ev}&X\ar@{=>}"a";"b"_F}$$
i.e.
$$-[f,F]=[f\nu,F(\S^{n-1}\nu)].$$
To see this, and in order to introduce the reader to ``track arguments'', we now prove by using diagrams that
$[f\nu,F(\S^{n-1}\nu)]$ is indeed the inverse of $[f,F]$.

By definition the sum $[f,F]+[f\nu,F(\S^{n-1}\nu)]$ is represented by diagram
\begin{equation*}\tag{a}
\xymatrix{S^n\ar[r]_<(.3){\S^{n-1}\mu}&S^n\vee S^n\ar[rr]_{\S^{n-1}(1\vee\nu)}&&S^n\vee
S^n\ar[rr]_{\S^{n-1}(f,f)}^{\;}="a"\ar@/^30pt/[rrr]^0_<(.285){\;}="b"&&
S^n_X\ar[r]_{ev}&X\ar@{=>}"a";"b"_{(F,F)}}
\end{equation*}
Since one-point unions of spaces are groupoid-enriched coproducts in $\C{Top}^*$ we have that the pasting of (a)
coincides with the pasting of
\begin{equation*}\tag{b}
\xymatrix{S^n\ar[r]_<(.3){\S^{n-1}\mu}&S^n\vee
S^n\ar[rr]_{\S^{n-1}(1,\nu)}&&
S^n\ar[r]_{\S^{n-1}f}^<(.6){\;}="a"\ar@/^30pt/[rr]^0_<(.285){\;}="b"&
S^n_X\ar[r]_{ev}&X\ar@{=>}"a";"b"_F}
\end{equation*}
The next step is crucial. We use the very important Remark \ref{*} to introduce a new track in diagram (b) (the
nil-track $N$ in (c) below) in such a way that the pasting of (b) and (c) is still the same.
\begin{equation*}\tag{c}
\xymatrix{S^n\ar[r]_<(.3){\S^{n-1}\mu}\ar@/^30pt/[rrr]^0_<(.68){\;}="d"&S^n\vee
S^n\ar[rr]_{\S^{n-1}(1,\nu)}^{\;}="c"&&
S^n\ar[r]_{\S^{n-1}f}^<(.6){\;}="a"\ar@/^30pt/[rr]^0_<(.285){\;}="b"&
S^n_X\ar[r]_{ev}&X\ar@{=>}"a";"b"_F\ar@{=>}"c";"d"^N}
\end{equation*}
This kind of argument will be very common throughout the whole paper. Indeed we use it here again in order to
remove the track $F$ from (c)  in such a way that the pasting of (c) is the same as the pasting of
\begin{equation*}\tag{d}
\xymatrix{S^n\ar[r]_<(.3){\S^{n-1}\mu}\ar@/^30pt/[rrr]^0_<(.69){\;}="d"&S^n\vee
S^n\ar[rr]_{\S^{n-1}(1,\nu)}^{\;}="c"&& S^n\ar[r]_{\S^{n-1}f}&
S^n_X\ar[r]_{ev}&X\ar@{=>}"c";"d"^N}
\end{equation*} 
The composite map in the bottom remains constant along diagrams (a)--(d). Moreover, we have also shown that the
pasting of these diagrams is also always the same. Therefore (a) and (d) represent the same element in
$\pi_{n,1}X$. We began with (a) which is $[f,F]+[f\nu,F(\S^{n-1}\nu)]$ and we finish with (d) which represents
the unit element $[0,0^\vc]$. Here we use the equivalence relation defining $\pi_{n,1}X$.

The reader can now easily check that $\partial$ is indeed a
homomorphism.
\end{proof}

\begin{defn}\label{cros}
A \emph{crossed module} $\partial\colon M\r N$ is a group homomorphism such that $N$ acts on the right of $M$ (the action will be denoted exponentially)
and the homomorphism
$\partial$ satisfies the following two properties $(m,m'\in M, n\in N)$:
\begin{enumerate}
\item $\partial(m^n)=-n+\partial(m)+n$,
\item $m^{\partial(m')}=-m'+m+m'$.
\end{enumerate}
A morphism $(f_0,f_1)\colon\partial\r\partial'$ between crossed modules $\partial\colon M\r N$ and
$\partial'\colon M'\r N'$ is a commutative square in the category of groups
$$\xymatrix{M\ar[d]_\partial\ar[r]^{f_1}&M'\ar[d]^{\partial'}\\N\ar[r]_{f_0}&N'}$$
such that for any $m\in M$ and $n\in N$ the formula
$f_1(m^n)=f_1(m)^{f_0(n)}$ holds. Such a morphism is a \emph{weak
equivalence} if it induces isomorphisms
$\ker\partial\cong\ker\partial'$ and
$\coker\partial\cong\coker\partial'$. The category of crossed
modules will be denoted by $\C{cross}$. A crossed module
$\partial\colon M\r N$ is \emph{$0$-free} if $N=\grupo{E}$ is a free
group.
\end{defn}

\begin{prop}\label{crosses}
The group $\pi_{1,0}X$ acts on $\pi_{1,1}X$ in such a way that
$\partial\colon\pi_{1,1}X\r\pi_{1,0}X$ is a crossed module.
Moreover, the induced map $\pi_{1,*}g$ in (\ref{inmor}) is a crossed module morphism.
\end{prop}

\begin{proof}
Let $\alpha\colon S^1\r S^1\vee S^1$ be any map inducing 
$\pi_1\alpha\colon\Z\r\grupo{a,b}\colon 1\mapsto -a+b+a$.
Any $x\in\pi_{1,0}X$ can
be identified with the homotopy class of a map $\tilde{x}\colon S^1\r S^1_X$. The automorphism
$$(-)^x\colon\pi_{1,1}X\To\pi_{1,1}X\colon [f,F]\mapsto[f,F]^x$$
is defined as follows: Let $[f,F]^x$ be given by the map
$$S^1\st{\alpha}\To S^1\vee S^1\st{(\tilde{x},f)}\To S^1_X$$
and the track
$$\xymatrix{&S^1\ar[r]^{\tilde{x}}_{\;}="d"&S^1_X\ar[rd]^{ev}&\\
S^1\ar@/^15pt/[ru]^0_{\;}="b"\ar[r]_-\alpha^<(.6){\;}="a"&S^1\vee S^1\ar[u]^{(1,0)}
\ar[r]_-{(\tilde{x},f)}^{\;}="c"&S^1_X\ar[r]_-{ev}&X\ar@{=>}"a";"b"^{N}\ar@{=>}"c";"d"_{(0^\vc,F)}}$$
Here $N$ is a nil-track. By using the elementary properties of nil-tracks in Theorem \ref{niltrackes} the reader can check that this is a well-defined action, independent
of the choice of $\alpha$. Equation (1) in Definition \ref{cros} is immediate. Let us now check that (2) holds. Consider $[f,F],[g,G]\in\pi_{1,1}X$.
On one hand $[f,F]^{\partial[g,G]}$ is represented by
\begin{equation*}\tag{a}
\xymatrix{&S^1\ar[r]^g_{\;}="d"&S^1_X\ar[rd]^{ev}&\\
S^1\ar@/^15pt/[ru]^0_{\;}="b"\ar[r]_-{\alpha}^<(.6){\;}="a"&S^1\vee S^1\ar[u]^{(1,0)}
\ar[r]_-{(g,f)}^{\;}="c"&S^1_X\ar[r]_-{ev}&X\ar@{=>}"a";"b"^N\ar@{=>}"c";"d"_{(0^\vc,F)}}
\end{equation*}
On the other hand $-[g,G]+[f,F]+[g,G]$ is represented by 
\begin{equation*}\tag{b}
\xymatrix{S^1\ar[r]_-{\alpha}&S^1\vee S^1\ar[r]_-{(g,f)}^<(.6){\;}="a"\ar@/^30pt/[rr]^0_{\;}="b"&S^1_X\ar[r]_-{ev}&X \ar@{=>}"a";"b"_{(G,F)}}
\end{equation*}
Factoring $(G,F)$ we can insert a $2$-cell to obtain
\begin{equation*}\tag{c}
\xymatrix{&S^1\ar[r]^g_{\;}="d"\ar@/^60pt/[rrd]^0_{\;}="b"&S^1_X\ar[rd]_{ev}^<(0){\;}="a"&\\
S^1\ar[r]_-{\alpha}&S^1\vee S^1\ar[u]^{(1,0)}
\ar[r]_-{(g,f)}^{\;}="c"&S^1_X\ar[r]_-{ev}&X \ar@{=>}"a";"b"^{G} \ar@{=>}"c";"d"_{(0^\vc,F)}}
\end{equation*}
so the pasting of (b) is the same as the pasting of (c).
Now we observe that the pasting of (a) and (c) coincide, hence we are done.
\end{proof}

\begin{defn}\label{qm}
A \emph{reduced quadratic module} $(\omega,\partial)$ is a sequence
of group homomorphisms $$N_\abb\otimes N_\abb\st{\omega}\To
M\st{\partial}\To N$$ such that, if $N\twoheadrightarrow
N_\abb\colon x\mapsto \set{x}$ is the projection onto the
abelianization, then the following equations hold for any $x,y\in N$
and $a,b\in M$,
\begin{enumerate}
\item $\partial\omega(\set{x}\otimes\set{y})=-x-y+x+y$,
\item $\omega(\set{\partial a}\otimes\set{\partial b})=-a-b+a+b$,
\item $\omega(\set{\partial a}\otimes\set{x}+\set{x}\otimes\set{\partial a})=0$.
\end{enumerate}
Moreover, it is a \emph{stable quadratic module} if the following condition, stronger than (3), is satisfied,
\begin{enumerate}\setcounter{enumi}{3}
\item $\omega(\set{x}\otimes\set{y}+\set{y}\otimes\set{x})=0$.
\end{enumerate}
Condition (4) says that $\omega\colon\otimes^2N_\abb\r M$ in a
stable quadratic module factors through the natural projection
$\bar{\sigma}\colon\otimes^2N_\abb\twoheadrightarrow
\hat{\otimes}^2N_\abb$. The factorization will also be denoted by
$\omega\colon\hat{\otimes}^2N_\abb\r M$ . We call
$(\omega,\partial\colon M\r N)$ \emph{$0$-free} if $N=\grupo{E}_\ni$
is a free nil-group.

A morphism of reduced or stable quadratic modules $(f_1,f_0)\colon(\omega,\partial)\r(\omega',\partial')$ is just a commutative diagram of the form
$$\xymatrix{N_\abb\otimes N_\abb\ar[r]^-\omega\ar[d]_{\otimes^2 (f_0)_\abb}&M\ar[r]^\partial\ar[d]^{f_1}&N\ar[d]^{f_0}\\
(N')_\abb\otimes(N')_\abb\ar[r]_-{\omega'}&M'\ar[r]_{\partial'}&N'}$$
It is a \emph{weak equivalence} if it induces isomorphisms
$\ker\partial\cong\ker\partial'$ and
$\coker\partial\cong\coker\partial'$. We will write $\C{rquad}$ and
$\C{squad}$ for the categories of reduced and symmetric quadratic
modules, respectively. See \cite{ch4c} IV.C.1.
\end{defn}

We now define a map $\omega$ as in Definition \ref{qm} for secondary homotopy groups $\pi_{n,*}X$ with $n\geq 2$. Consider the diagram
$$\xymatrix{S^n\ar[rr]_{\S^{n-1}\beta}^{\;}="a"\ar@/^30pt/[rr]^0_{\;}="b"&&S^n\vee S^n\ar@{=>}"a";"b"_{B}}$$
where $\beta\colon S^1\r S^1\vee S^1$ is given such that $(\pi_1\beta)_\ni(1)=-a-b+a+b\in\grupo{a,b}_\ni$ is the commutator. The track $B$ is determined by
$\hopf(B)=-a\otimes b \in\otimes^2\Z[a,b]$. 
Given $x\otimes y\in\otimes^2_n(\pi_{n,0}X)_\abb$ let
$\tilde{x},\tilde{y}\colon S^1\r\vee_{\L^nX}S^1$ be maps with
$(\pi_1\tilde{x})_\abb(1)=x$ and $(\pi_1\tilde{y})_\abb(1)=y$. Then
the diagram
\begin{equation*}
\xymatrix{S^n\ar[rr]_{\S^{n-1}\beta}^{\;}="a"\ar@/^30pt/[rr]^0_{\;}="b"&&S^n\vee
S^n\ar@{=>}"a";"b"_{B}\ar[rr]_{\S^{n-1}(\tilde{x},\tilde{y})}&&S^n_X\ar[r]_{ev}&X}
\end{equation*}
represents an element
\begin{equation}\label{omede}
\omega(x\otimes y)=[(\tilde{x},\tilde{y})\beta,ev\,(\S^{n-1}(\tilde{x},\tilde{y}))B]\in\pi_{n,1}X.
\end{equation}

\begin{prop}\label{redquades}
For $n\geq 2$ the homomorphism of groups 
$$\omega\colon(\pi_{n,0}X)_\abb\otimes(\pi_{n,0}X)_\abb\To\pi_{n,1}X$$
given by (\ref{omede}) is well defined. Moreover,
 $(\omega,\partial)$ is a reduced quadratic module for $n=2$ and a stable quadratic module for $n\geq 3$.
Furthermore, (\ref{inmor}) is a reduced quadratic module
homomorphism for all $n\geq 2$.
\end{prop}

\begin{proof}
Given $x\in\otimes^2_n(\pi_{n,0}X)_\abb$ we alternatively define
$$\omega(x)=[\omega(x)_1,ev\,\omega(x)_2]\in\pi_{n,1}X$$
by choosing $\omega(x)_1$ and
$\omega(x)_2\colon(\omega(x)_1)_{ev}\rr 0$ as in diagram
$$\xymatrix{S^n\ar[rr]_{\S^{n-1}\omega(x)_1}^{\;}="a"\ar@/^25pt/[rr]^0_{\;}="b"&&S^n_X\ar[r]_{ev}&X\ar@{=>}"a";"b"^{\omega(x)_2}}$$
Here $\omega(x)_1\colon S^1\r\vee_{\L^nX}S^1$ is a map with
$(\pi_1\omega(x)_1)_\ni(1)=\partial(x)$ and $$\omega(x)_2\colon
\S^{n-1}\omega(x)_1\rr 0$$ is the unique track with
$\hopf(\omega(x)_2)=-x$. Such a track exists and is unique by
Theorem \ref{propi}. The elementary properties of nil-tracks and
more generally of the Hopf invariant for tracks, see Theorem
\ref{propi}, show that the element $\omega(x)$ is indeed
well-defined and this definition of $\omega$ coincides with the
definition given by (\ref{omede}). Axiom (1) in Definition \ref{qm}
is automatically satisfied. The bilinearity of $\omega$ follows
from Theorem \ref{propi} (1) and the following claim:
\begin{itemize}
\item[(*)] Given $[f,F]\in\pi_{1,n}X$ and $x\in\otimes^2_n(\pi_{n,0}X)_\abb$ then the sum $[f,F]+\omega(x)=[g,G]\in\pi_{n,1}X$ is
represented by a map $g$ with $(\pi_1 g)_\ni(1)=\partial[f,F]+\partial(x)$ and $G=F\vc(ev\,\bar{G})$ where $\bar{G}\colon\S^{n-1}g\rr\S^{n-1}f$
is the unique track with $\hopf(\bar{G})=-x$.
\end{itemize}
Indeed $[f,F]+\omega(x)=[(f,\omega(x)_1)\mu,(F,ev\,\omega(x)_2)(\S^{n-1}\mu)]$ and we can suppose $g=(f,\omega(x)_1)\mu$. 
The pasting of the following diagram
is trivial
$$\xymatrix{S^n\ar@/_5pt/@{=}[rd]^{\;}="d"\ar[r]^(.4){\Sigma^{n-1}\mu}&S^n\vee S^n\ar[d]|{(1,0)}_<(.05){\;}="c"
\ar[rrr]|<(.45){(\S^{n-1}f,\S^{n-1}\omega(x)_1)}^<(.6){\;}="a"_{\;}="h"\ar@/^40pt/[rrrr]^0_{\;}="b"&&&
S^n_X\ar[r]_{ev}&X\ar@{=>}"a";"b"_{(F,ev\,\omega(x)_2)}\\
&S^n\ar@/_10pt/[rrru]|{\S^{n-1}f}_<(.6){\;}="f"^{\;}="g"\ar@/_25pt/[rrrru]_0^{\;}="e"\ar@{<=}"c";"d"_N\ar@{=>}"e";"f"_{F^\vi}\ar@{=>}"g";"h"^{(0^\vc,
\omega(x)_2^\vi)}}$$
since, provided we remove the nil-track $N$ from the diagram, the pasting of the two lower tracks  is the inverse
of the upper track, and by Remark \ref{*} the nil-track $N$ can be freely removed from this diagram without
altering the result of the pasting.
Moreover, it is not difficult to see by using Theorem \ref{propi} that
$$\hopf(((\S^{n-1}f)N^\vi)\vc((0^\vc,
\omega(x)_2)(\S^{n-1}\mu)))=-x.$$ therefore
$\bar{G}=((\S^{n-1}f)N^\vi)\vc((0^\vc, \omega(x)_2)(\S^{n-1}\mu))$
and the claim follows.


The map $\beta\colon S^1\r S^1\vee S^1$ induces the commutator of the generators on $\pi_1$, therefore by
standard ``universal example" arguments, by the elementary properties of nil-tracks, and by Remark \ref{*}, we see
that the commutator
$-[f,F]-[g,G]+[f,F]+[g,G]\in\pi_{n,1}X$ is represented by the following
diagram
$$\xymatrix{S^n\ar[rr]_{\S^{n-1}\beta}&&S^n\vee S^n\ar@/^30pt/[rrr]^0_<(.28){\;}="a"\ar[rr]_{\S^{n-1}(f,g)}^{\;}="b"\ar@{<=}"a";"b"^{(F,G)}&&S^n_X\ar[r]_{ev}&X}$$
By using Remark \ref{*} twice, as in the proof of Proposition \ref{grupoes}, we deduce that the pasting of
the previous diagram coincides with the pasting of
\begin{equation*}
\xymatrix{S^n\ar[rr]_{\S^{n-1}\beta}^{\;}="a"\ar@/^30pt/[rr]^0_{\;}="b"&&S^n\vee
S^n\ar@{=>}"a";"b"_{B}\ar[rr]_{\S^{n-1}(f,g)}&&S^n_X\ar[r]_{ev}&X}
\end{equation*}
Therefore (2) in Definition \ref{qm} is satisfied.

If $n\geq 3$ equation (3) in Definition \ref{qm} is an immediate
consequence of the fact that
$\bar{\sigma}(\set{x}\otimes\set{y}+\set{y}\otimes\set{x})=0$. Let
us now check (3) in Definition \ref{qm} in case $n=2$. Suppose that
we have $[f,F]\in\pi_{2,1}X$ and $x\in\pi_{2,0}X$. We choose
$\tilde{x}\colon S^1\r\vee_{\L^2X}S^1$ such that
$(\pi_1\tilde{x})_\ni(1)=x$. Then by claim (*) and (\ref{omede}) we have that
$\omega(\set{\partial[f,F]}\otimes\set{x}+\set{x}\otimes\set{\partial[f,F]})$
is represented by the diagram
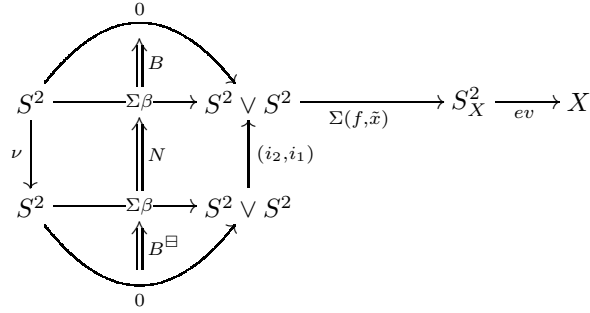
\begin{equation*}
\xymatrix{S^2\ar[d]_\nu\ar[rr]|{\S\beta}^{\;}="a"_{\;}="d"\ar@/^30pt/[rr]^0_{\;}="b"&&S^2\vee
S^2\ar@{=>}"a";"b"_{B}\ar[rr]_{\S(f,\tilde{x})}&&S^2_X\ar[r]_{ev}&X\\
S^2\ar[rr]|{\S\beta}^{\;}="c"_{\;}="f"\ar@/_30pt/[rr]_0^{\;}="e"&&S^2\vee
S^2\ar[u]_{(i_2,i_1)}&\ar@{=>}"c";"d"_N\ar@{=>}"e";"f"_{B^\vi}}
\end{equation*}
Here $N$ is a nil-track. By Remark \ref{*} the pasting of this diagram coincides with the pasting of 
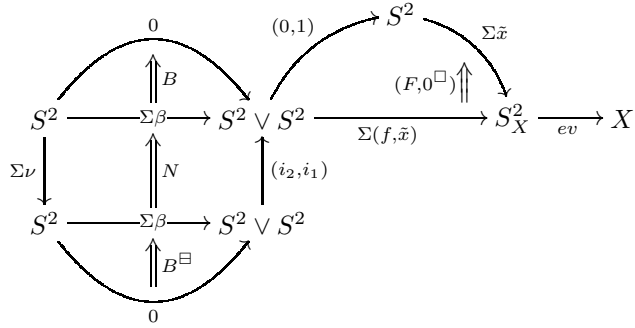
\begin{equation*}
\xymatrix{&&&S^2\ar@/^10pt/[rd]^{\S\tilde{x}}_{\;}="h"\\S^2\ar[d]_{\S\nu}\ar[rr]|{\S\beta}^{\;}="a"_{\;}="d"\ar@/^30pt/[rr]^0_{\;}="b"&&S^2\vee
S^2\ar@/^10pt/[ru]^{(0,1)}\ar@{=>}"a";"b"_{B}\ar[rr]_{\S(f,\tilde{x})}^<(.76){\;}="g"&&S^2_X\ar[r]_{ev}&X\\
S^2\ar[rr]|{\S\beta}^{\;}="c"_{\;}="f"\ar@/_30pt/[rr]_0^{\;}="e"&&S^2\vee
S^2\ar[u]_{(i_2,i_1)}&\ar@{=>}"c";"d"_N\ar@{=>}"e";"f"_{B^\vi}\ar@{=>}"g";"h"^{(F,0^\vc)}}
\end{equation*}
and again by Remark \ref{*} the pasting of this diagram is the same as the pasting of
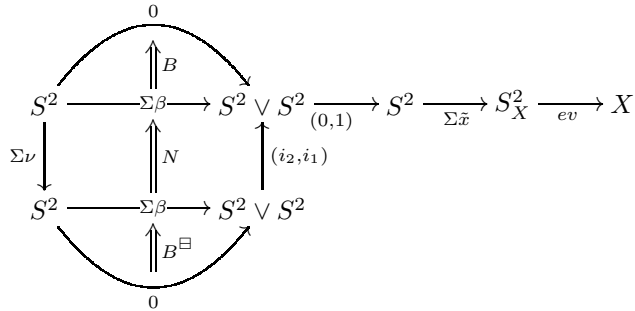
\begin{equation*}
\xymatrix{S^2\ar[d]_{\S\nu}\ar[rr]|{\S\beta}^{\;}="a"_{\;}="d"\ar@/^30pt/[rr]^0_{\;}="b"&&S^2\vee
S^2\ar@{=>}"a";"b"_{B}\ar[r]_{(0,1)}&S^2\ar[r]_{\S\tilde{x}}^<(.76){\;}="g"&S^2_X\ar[r]_{ev}&X\\
S^2\ar[rr]|{\S\beta}^{\;}="c"_{\;}="f"\ar@/_30pt/[rr]_0^{\;}="e"&&S^2\vee
S^2\ar[u]_{(i_2,i_1)}&\ar@{=>}"c";"d"_N\ar@{=>}"e";"f"_{B^\vi}}
\end{equation*}
By using Theorem \ref{propi} one can readily check that $(0,1)(B\vc
N\vc((i_2,i_1)B^\vi(\S\nu)))$ is a nil-track, and therefore this
diagram represents the trivial element in $\pi_{2,1}X$.
\end{proof}

For $n\geq 0$ we define the category $\C{cross}(n)$ as follows. 
\begin{equation}\label{crossn}
\C{cross}(n)=\left\{\begin{array}{ll}\C{grd}^*,&\text{pointed groupoids if $n=0$};\\&\\
\C{cross},&\text{crossed modules if $n=1$};\\&\\\C{rquad},&\text{reduced quadratic modules if $n=2$};\\&\\
\C{squad},&\text{stable quadratic modules if $n\geq 3$}.\end{array}\right.
\end{equation}

\begin{thm}
Secondary homotopy groups are well-defined functors
$$\pi_{n,*}\colon\C{Top}^*\To\C{cross}(n),\;\;\;n\geq 0.$$
\end{thm}

This result generalizes the well-known fact on classical homotopy
groups which are functors
$$\pi_n\colon\C{Top}^*\To\C{group}(n),\;\;n\geq0,$$
where
\begin{equation}\label{groupn}
\C{group}(n)=\left\{
                 \begin{array}{ll}
                   \C{Set}^*, & \hbox{pointed sets if $n=0$;} \\
&\\
                   \C{Gr}, & \hbox{groups if $n=1$;} \\
&\\
                   \C{Ab}, & \hbox{abelian groups if $n\geq 2$.}
                 \end{array}
               \right.
\end{equation}

Moreover, we have functors $(n\geq 0)$
\begin{equation}\label{hi}
\begin{array}{c}
h_0\colon\C{cross}(n)\To\C{group}(n),\\
h_1\colon\C{cross}(n)\To\C{group}(n+1).
\end{array}
\end{equation}
The functor $h_0$ is defined as
the cokernel of the group homomorphism $\partial\colon M\r N$ for a
crossed module $\partial$ or a reduced or stable quadratic module
$(\omega,\partial)$, and for $\C{G}$ a pointed groupoid
$h_0\C{G}=\iso(\C{G})$ is the pointed set of isomorphism classes of
objects. Similarly $h_1$ is the kernel of $\partial\colon M\r N$ for
crossed modules and reduced and stable quadratic modules, and
$h_1\C{G}=\aut_\C{G}(*)$ is the automorphism group of the distinguished
object. 
In particular a morphism $f$ in $\C{cross}(n)$ is a
\emph{weak equivalence} for $n\geq 1$ if and only if $h_0f$ and $h_1f$ are
isomorphisms.

In Proposition \ref{exacta} below we show that there are natural isomorphisms $(n\geq 0)$ 
\begin{equation}\label{full1}
h_0\pi_{n,*}X\cong\pi_nX\text{ and } h_1\pi_{n,*}X\cong \pi_{n+1}X.
\end{equation}

Our definition of $\pi_{n,*}X$ above is a ``singular'' and hence
functorial version of secondary homotopy groups. For many purposes
it suffices to consider smaller models of $\pi_{n,*}X$ by choosing a
subset of $\L^nX$ which generates $\pi_nX$. Let
us make precise this observation.

\begin{prop}\label{menos}
Let $X$ be a pointed space. If $E_0\r X$ is a pointed map
between pointed sets then there is a unique pointed groupoid
$\pi_{0,*}(X,E_0)$ with object set $E_0$ endowed with a full and
faithful functor
$$\pi_{0,*}(X,E_0)\To\pi_{0,*}X$$
given by $E_0\r X$ on object sets. This morphism of pointed
groupoids is a weak equivalence provided any component of $X$ has
points in the image of $E_0$. Moreover, a map of pointed sets
$E_1\r\L X$ induces a crossed module morphism by the pull-back
$$\xymatrix@C=20pt{\pi_{1,1}(X,E_1)\ar[d]_\partial\ar[r]\ar@{}[rd]|{\text{pull}}&\pi_{1,1}X\ar[d]^\partial\\
\;\;\;\;\;\;\;\;\;\;\;\pi_{1,0}(X,E_1)=\grupo{E_1}\ar[r]&\grupo{\L
X}=\pi_{1,0}X\;\;\;\;\;\;\;\;\;\;\;\;\;\;}$$ which is a weak
equivalence 
$$\pi_{1,*}(X,E_1)\st{\sim}\To\pi_{1,*}X$$
provided the loops in the image of $E_1$ generate the group
$\pi_1X$. Furthermore, for $n\geq 2$ the a map of pointed sets
$E_n\r\L^n X$ induces a reduced (stable if $n\geq 3$) quadratic module morphism by the pull-back
$$\xymatrix@C=30pt{\;\;\;\;\;\;\;\;\;\;\;\;\;\;\;\;\;\otimes^2(\pi_{n,0}(X,E_n))_\abb=\otimes^2\Z[E_n]\ar[d]_\omega\ar[r]
\ar@/_40pt/[dd]_\partial&
\otimes^2\Z[\L^nX]=\otimes^2(\pi_{n,0}X)_\abb\ar[d]^\omega\;\;\;\;\;\;\;\;\;\;\;\;\;\;\;\;\;\;\\\pi_{n,1}(X,E_n)\ar[d]_\partial\ar[r]\ar@{}[rd]|{\text{pull}}&\pi_{n,1}X\ar[d]^\partial\\
\;\;\;\;\;\;\;\;\;\;\;\pi_{n,0}(X,E_n)=\grupo{E_n}_\ni\ar[r]&\grupo{\L^n
X}_\ni=\pi_{n,0}X\;\;\;\;\;\;\;\;\;\;\;\;\;\;}$$ which is a weak
equivalence 
$$\pi_{n,*}(X,E_n)\st{\sim}\To\pi_{n,*}X,\;\;n\geq 2$$
provided the $n$-loops in the image of $E_n$ generate the abelian
group $\pi_nX$.
\end{prop}

This proposition can be used to reduce the number of generators of a
secondary homotopy group, as one can check in the following
example.

\begin{rem}
So far we have not computed any secondary homotopy group. Now, with
the help of Proposition \ref{menos} we give a small model for the
secondary homotopy group $\pi_{n,*}(\vee_ES^n)$ of a wedge of
spheres indexed by the pointed set $E$. For this we notice that
there is a pointed inclusion $E\subset\L^n(\vee_ES^n)$ sending $e\in
E-\set{*}$ to the inclusion of the corresponding factor of the wedge
$S^n\subset\vee_ES^n$. Then we have a weak
equivalence
$$\pi_{n,*}(\vee_ES^n,E)\st{\sim}\To\pi_{n,*}(\vee_ES^n),\;\;n\geq 1.$$
For $n=1$ one easily checks that $\pi_{1,*}(\vee_ES^1,E)$ is 
$$\pi_{1,1}(\vee_ES^1,E)=0\st{\partial}\To\pi_{1,0}(\vee_ES^1,E)=\grupo{E}.$$
For $n=2$ the reduced quadratic module $\pi_{2,*}(\vee_ES^2,E)$ is given by the following diagram, see
(\ref{exagam})
$$\xymatrix{\otimes^2(\pi_{n,0}(\vee_ES^n,E))_\abb\ar[r]^\omega\ar@{=}[d]&\pi_{n,1}(\vee_ES^n,E)\ar@{=}[d]\ar[r]^\partial&
\pi_{n,0}(\vee_ES^n,E)\ar@{=}[d]\\\otimes^2\Z[E]\ar@{=}[r]&\otimes^2\Z[E]\ar[r]^{\partial}&\grupo{E}_\ni}$$
This follows from the fact that the next diagram is a pull-back
$$\xymatrix@C=30pt{\otimes^2\Z[E]\ar@{=}[d]\ar@{^{(}->}[r]
\ar@/_40pt/[dd]&
\otimes^2\Z[\L^2\vee_ES^2]=\otimes^2(\pi_{2,0}\vee_ES^2)_\abb\ar[d]^\omega\;\;\;\;\;\;\;\;\;\;\;\;\;\;\;\;\;\;\;\;\;\;\;\;\;\;\;\;\\
\otimes^2\Z[E]\ar[d]_\partial\ar@{^{(}->}[r]^\phi\ar@{}[rd]|{\text{pull}}&\pi_{2,1}\vee_ES^2\ar[d]^\partial\\
\grupo{E}_\ni\ar@{^{(}->}[r]&\grupo{\L^2
\vee_ES^2}_\ni=\pi_{2,0}\vee_ES^2\;\;\;\;\;\;\;\;\;\;\;\;\;\;\;\;\;\;\;\;\;\;\;\;\;\;}$$
Here the homomorphism $\phi$ is defined as follows. Given
$x\in\otimes^2\Z[E]$ the element $\phi(x)=[\phi_1(x),\phi_2(x)]$ is
given by a map $$\phi_1(x)\colon S^1\st{\bar{\phi}(x)}\To
\vee_ES^1\subset\vee_{\L^2\vee_ES^2}S^1$$ with
$(\pi_1\bar{\phi}(x))_\ni(1)=\partial(x)$ and the unique track
$$\xymatrix{S^2\ar[r]_{\S \phi_1(x)}^<(.98){\;\;\;\;\;}="a"\ar@/^25pt/[rr]^0_{}="b"&S^2_{\vee_ES^2}\ar[r]_{ev}&\vee_ES^2\ar@{=>}"a";"b"_{\phi_2(x)}}$$
with Hopf invariant $\hopf(\phi_2(x))=-x$. Here we use the fact that the
composite $ev(\S\phi_1(x))=\S\bar{\phi}(x)$ is a suspension. For
$n\geq 3$ the stable quadratic module $\pi_{n,*}(\vee_ES^n,E)$ is
given by the following diagram, see (\ref{exagam}).
$$\xymatrix{\otimes^2(\pi_{n,0}(\vee_ES^n,E))_\abb\ar[r]^-\omega\ar@{=}[d]&\pi_{n,1}(\vee_ES^n,E)\ar@{=}[d]\ar[r]^\partial&
\pi_{n,0}(\vee_ES^n,E)\ar@{=}[d]\\\otimes^2\Z[E]\ar[r]^{\bar{\sigma}}&\hat{\otimes}^2\Z[E]\ar[r]^{\partial}&\grupo{E}_\ni}$$
This can be easily checked as in the case $n=2$ by using the Hopf
invariant for tracks.
\end{rem}

\section{Homotopy groups of fibers}

We first obtain by secondary homotopy groups the classical homotopy
groups $\pi_nX$ as in the next result.

\begin{prop}\label{exacta}
For all $n\geq 1$ there is a natural exact sequence of groups
$$\pi_{n+1}X\st{\iota}\hookrightarrow\pi_{n,1}X\st{\partial}\To\pi_{n,0}X\st{ q }\twoheadrightarrow\pi_nX,$$
where $ q $ sends a basis element of $\pi_{n,0}X$, which is a map $f\colon S^n\r X$, to its homotopy class in $\pi_nX$; and $\iota$ carries the
homotopy class of $f\colon S^{n+1}\r X$ to the element $[0,pf]\in\pi_{n,1}X$, where 
$p\colon IS^n\r\S S^n=S^{n+1}$ is the obvious projection. 
\end{prop}


\begin{proof}
Obviously $ q $ is surjective. Any element $x\in\pi_{n,0}X$ is
represented by a map $\tilde{x}\colon S^1\r\vee_{\L^nX}S^1$, i.e.
$(\pi_1\tilde{x})_\ni(1)=x$. It is immediate to notice that $ q (x)$
is the homotopy class of $\tilde{x}_{ev}\colon S^n\r X$. If $ q
(x)=0$ then there exists a track $H\colon\tilde{x}_{ev}\rr 0$, and
the pair $[\tilde{x},H]\in\pi_{n,1}X$ satisfies
$\partial[\tilde{x},H]- q (x)$. It is immediate to notice that $ q
\partial=0$ and $\partial\iota=0$. The injectivity of $\iota$ is
also easy to check, actually $\pi_{n+1}X$ is isomorphic to the
subgroup of $\pi_{n,1}X$ given by the elements which can be
represented with a $0$ in the first coordinate. Finally suppose that
for some $[f,F]\in\pi_{n,1}X$ we have $\partial[f,F]=0$, then the
nil-track $N\colon0\rr\S^{n-1}f$ is defined and $[f,F]=[0,F\vc N]$,
hence we are done.
\end{proof}


We now introduce the (algebraic) fiber of a map in $\C{cross}(n)$ for $n\geq 1$.

\begin{defn}\label{fibers}
Let $f\colon\partial\r\partial'$ be a crossed module morphism
\begin{equation*}
\xymatrix{M\ar[d]_\partial\ar[r]^{f^1}&M'\ar[d]^{\partial'}\\N\ar[r]_{f^0}&N'}
\end{equation*}
We define the \emph{fiber $\fib(f)$} as the crossed module
$\fib(f)\colon\fib_1(f)\r\fib_0(f)$ where $\fib_0(f)$ is the
following pull-back
\begin{equation*}
\xymatrix{\fib_0(f)\ar[d]_{\bar{\partial}'}\ar[r]^{\bar{f}^0}\ar@{}[rd]|{\text{pull}}&M'\ar[d]^{\partial'}\\N\ar[r]_{f^0}&N'}
\end{equation*}
$\fib_1(f)=M$ and the homomorphism $\fib(f)\colon\fib_1(f)\r\fib_0(f)$ is induced by $(\partial,f^1)\colon M\r N\times M'$.
The action of $\fib_0(f)$ on $\fib_1(f)$ is the pull-back along $\bar{\partial}'$ of the action of $N$ on $M$. The axioms of a crossed module are easily verified.
There is also a natural crossed module morphism $\jmath\colon\fib(f)\r\partial$ given by the square
\begin{equation*}
\xymatrix{\fib_1(f)\ar[d]_{\fib(f)}\ar@{=}[r]&M\ar[d]^{\partial}\\\fib_0(f)\ar[r]_{\bar{\partial}'}&N}
\end{equation*}

Let $f\colon(\omega,\partial)\r(\omega',\partial')$ be now a reduced/stable quadratic module morphism
\begin{equation*}
\xymatrix{\otimes^2N_\abb\ar[d]^\omega\ar[r]^{\otimes^2f^0_\abb}&\otimes^2(N')_\abb\ar[d]_{\omega'}\\M\ar[d]_\partial\ar[r]^{f^1}&M'\ar[d]^{\partial'}\\N\ar[r]_{f^0}&N'}
\end{equation*}
The \emph{fiber $\fib(f)$} is a reduced/stable quadratic module
$$\otimes^2(\fib_0(f))_\abb\To\fib_1(f)\st{\fib(f)}\To\fib_0(f)$$
where $\fib(f)\colon\fib_1(f)\To\fib_0(f)$ is defined as in the
crossed module case and the first homomorphism is the composite
$$\otimes^2(\fib_0(f))_\abb\st{\otimes^2\bar{\partial}'_\abb}\To\otimes^2N_\abb\st{\omega}\To M.$$
The natural reduced/stable quadratic module morphism $\jmath\colon\fib(f)\r(\omega,\partial)$
is also defined as above.
\end{defn}

\begin{lem}\label{ex}
Let $f\colon\partial\r\partial'$ be a morphism of crossed modules, then there is an exact sequence
$$h_1\fib(f)\st{h_1\jmath}\hookrightarrow h_1\partial\st{h_1f}\To h_1\partial'\st{\delta}\To
h_0\fib(f)\st{h_0\jmath}\To h_0\partial\st{h_0f}\To h_0\partial'.$$
This exact sequence is natural in $f$. Moreover, it is also available for reduced or stable quadratic module 
morphisms $f\colon(\omega,\partial)\r(\omega',\partial')$.
\end{lem}

\begin{proof}
The homomorphism $\delta$ is determined
by the inclusion $M'\hookrightarrow N\times M'\colon m'\mapsto (0,m')$. The proof of the exactness is a simple exercise.
\end{proof}

\begin{thm}
Let $f\colon X\r Y$ be a map between pointed spaces and let $F_f$ be the homotopy fiber of $f$. Then for all $n\geq 1$ there is a natural morphism in $\C{cross}(n)$
$$\xi\colon\pi_{n,*}F_f\To\fib(\pi_{n,*}f)$$
which induces an isomorphism
\begin{equation*}\tag{1}
\pi_nF_f\cong h_0\fib(\pi_{n,*}f)
\end{equation*}
and an exact sequence
\begin{equation*}\tag{2}
\pi_{n+2}Y\To\pi_{n+1}F_f\twoheadrightarrow h_1\fib(\pi_{n,*}f),
\end{equation*}
where the first arrow is the boundary homomorphism in the long exact
sequence in homotopy. By using the isomorphism $(1)$ above and
Proposition \ref{exacta} we can naturally identify the exact
sequence in Lemma \ref{ex} extended on the left by the exact
sequence $(2)$ with the following piece of the long exact sequence
of homotopy groups
$$\pi_{n+2}Y\r\pi_{n+1}F_f\r\pi_{n+1}X\r\pi_{n+1}Y\r\pi_nF_f\r\pi_nX\r\pi_nY.$$
\end{thm}

\begin{proof}
Recall that $F_f$ is a pull-back
\begin{equation*}
\xymatrix{F_f\ar[d]_{\bar{e}}\ar[r]^{\bar{f}}\ar@{}[rd]|{\text{pull}}&Y^I\ar[d]^{ev_0}\\X\ar[r]_f&Y}
\end{equation*}
where $Y^I$ is the space of based maps $([0,1],1)\r(Y,*)$ and $ev_0$ is the evaluation at $0\in [0,1]$.

The morphism $\xi$ consists of two morphisms, the upper one is $$\xi^1=\pi_{n,1}\bar{e}\colon\pi_{n,1}F_f\r\pi_{n,1}X=\fib_1(\pi_{n,*}f).$$

We now construct the map
$\xi_0\colon\pi_{n,0}F_f\r\fib_0(\pi_{n,*}f)$. For this we consider
on the one hand the morphism
$\pi_{n,0}\bar{e}\colon\pi_{n,0}F_f\r\pi_{n,0}X$ induced by $f$. On
the other hand we define a homomorphism
$$\bar{\xi}\colon\pi_{n,0}F_f\r\pi_{n,1}Y$$ as follows: an element $z\in\pi_{n,0}F_f$
is represented by a map $\tilde{z}\colon S^1\r\vee_{\L^nF_f}S^1$ with $(\pi_1\tilde{z})(1)=z$ if $n=1$ or
$(\pi_1\tilde{z})_\ni(1)=z$ if $n\geq 2$. The map $$S^n\st{\S^{n-1}\tilde{z}}\To S^n_{F_f}\st{\S^n\L^n\bar{f}}\To S^n_{Y^I}\st{ev}\To Y^I$$
has an adjoint $$ad(ev(\S^n\L^n\bar{f})(\S^{n-1}\tilde{z}))\colon IS^n\To Y,$$
this adjoint represents a track $ad(ev(\S^n\L^n\bar{f})(\S^{n-1}\tilde{z}))\colon((\S\L^n(f\bar{e}))\tilde{z})_{ev}\rr0$, and
$$\bar{\xi}(z)=[(\S\L^n(f\bar{e}))\tilde{z},ad(ev(\S^n\L^n\bar{f})(\S^{n-1}\tilde{z}))]\in\pi_{n,1}Y.$$
It is immediate to check that $\pi_{n,0}\bar{e}$ and $\bar{\xi}$ define a homomorphism to the pull-back
$$\xi^0=(\pi_{n,0}\bar{e},\bar{\xi})\colon\pi_{n,0}F_f\To\fib_0(\pi_{n,*}f).$$
Now it is easy to check that $\xi$ is indeed a morphism in $\C{cross}(n)$.

By Proposition \ref{exacta} and Lemma \ref{ex} we obtain from $\xi$ a diagram with exact rows
\begin{scriptsize}
$$\xymatrix@C=10pt{\pi_{n+2}Y\ar[r]&\pi_{n+1}F_f\ar[r]\ar[d]^{\xi_*}&\pi_{n+1}X\ar[r]\ar[d]^\cong&\pi_{n+1}Y\ar[d]^\cong\ar[r]&\pi_nF_f\ar[d]^{\xi_*}\ar[r]&\pi_nX
\ar[d]^\cong\ar[r]&\pi_nY\ar[d]^\cong\\
&h_1\fib(\pi_{n,*}f)\ar@{^{(}->}[r]&
h_1\pi_{n,*}X\ar[r]&h_1\pi_{n,*}Y\ar[r]^\delta&h_0\fib(\pi_{n,*}f)\ar[r]&h_0\pi_{n,*}X\ar[r]&h_0\pi_{n,*}Y}$$
\end{scriptsize}
It is easy to see that this diagram commutes, and hence the theorem follows from the five lemma.
\end{proof}

\begin{cor}
Let $f\colon X\r Y$ be a map between pointed spaces and let $F_f$ be the homotopy fiber of $f$. If $\pi_{n+2}f\colon\pi_{n+2}X\twoheadrightarrow\pi_{n+2}Y$
is surjective then there is a weak equivalence in $\C{cross}(n)$, $n\geq 1$,
$$\xi\colon\pi_{n,*}F_f\st{\sim}\To\fib(\pi_{n,*}f).$$
\end{cor}

\section{Suspension and loop functors}

Homotopy groups $\pi_nX$ are objects in the category $\C{group}(n)$,
$n\geq 0$, see (\ref{groupn}). There are forgetful functors
\begin{equation}\label{loops1}
\phi_n\colon\C{group}(n)\To\C{group}(n-1)
\end{equation} given by $\phi_n=1_\C{Ab}$ for $n\geq 3$ and by the
obvious forgetful functors $\phi_2\colon\C{Ab}\To\C{Gr},$
$\phi_1\colon\C{Gr}\To\C{Set}^*.$

It is a classical result that for any pointed space $X$ there are
natural isomorphisms $n\geq 0$ $$\L\colon\pi_n\L
X\cong\phi_{n+1}\pi_{n+1}X\text{ in }\C{group}(n).$$ The analogue of
this isomorphism for secondary homotopy groups is as follows.

There are forgetful functors
\begin{equation}\label{loops2}
\phi_n\colon\C{cross}(n)\To\C{cross}(n-1),
\end{equation}
see (\ref{crossn}),
given by $\phi_n=1_\C{squad}$ for $n\geq 4$ and by the functors
\begin{equation}\label{loops}
\begin{array}{c}
\phi_3\colon\C{squad}\To\C{rquad},\\
\phi_2\colon\C{rquad}\To\C{cross},\\
\phi_1\colon\C{cross}\To\C{grd}^*.
\end{array}
\end{equation}
The functor $\phi_3$ in (\ref{loops}) is obvious, since stable
quadratic modules are special reduced quadratic modules. Given a
reduced quadratic module $(\omega,\partial)$ we have
$\phi_2(\omega,\partial)=\partial\colon M\r N$ in $\C{cross}$, with
the action of $N$ on $M$ defined by
$$m^n=m+\omega(\set{\partial m}\otimes \set{n}).$$ Finally if
$\partial\colon M\r N$ is a crossed module then the pointed groupoid $\phi_1\partial$ in $\C{grd}^*$
has $N$ as a set of objects. 
Moreover the set of all morphisms in $\phi_1\partial$ is the
semidirect product $N\ltimes M$, which is the group structure on the set $N\times M$ defined by the formula
$$(n,m)+(n',m')=(n+n',m^{n'}+m'),$$
 and the structure maps of the groupoid (identities, source and target)
$$N\st{i}\r N\ltimes M\mathop{\rightrightarrows}\limits_t^s N$$ 
are $i(n)=(n,0)$, $s(n,m)=n$ and $t(n,m)=n+\partial m$. The composition law $\circ$ is determined by the formula
$$(n+\partial m, m')\circ(n,m)=(n,m+m').$$

The forgetful functors $\phi_n$ in (\ref{loops2}) clearly
commute with 
$h_0$ and $h_1$ in (\ref{hi}), that is,
\begin{equation}\label{conmuta1}
h_i\phi_n=\phi_n h_i,\;\;n\geq1,i=0,1.
\end{equation}

\begin{thm}\label{adiso}\label{adeq}
There is a natural weak equivalence in $\C{cross}(n)$
$$\L\colon\pi_{n,*}\L X\To\phi_{n+1}\pi_{n+1,*}X,\;\;n\geq0,$$ which induces
the isomorphism $\L\colon\pi_n\L X\cong\phi_{n+1}\pi_{n+1}X$ in
$h_0$ and $-\L\colon\pi_{n+1}\L X\cong\phi_{n+2}\pi_{n+2}X$ in
$h_1$. This weak equivalence is an isomorphism for $n\geq 2$.
\end{thm}

\begin{proof}
Let us first consider the case $n\geq 3$.

We have $\pi_{n,0}\L X=\grupo{\L^{n+1}X}_\ni=\pi_{n+1,0}X$. We
define a group homomorphism $\pi_{n,1}\L X\r\pi_{n+1,1}X$ sending
$[f,F]$ with $f\colon S^1\r\vee_{\L^{n+1}X}S^1$ and
$$\xymatrix{S^n\ar[r]_{\S^{n-1}f}^<(.98){\;\;\;\;\;}="a"\ar@/^25pt/[rr]^0_{}="b"&S^n_{\L X}
\ar[r]_{ev}&\L X\ar@{=>}"a";"b"_F}$$ to $[f,ad(F)]$ where $ad(F)$ is
the adjoint track
$$\xymatrix{S^{n+1}\ar[r]_{\S^n f}^<(.98){\;\;\;\;\;}="a"\ar@/^25pt/[rr]^0_{}="b"&S^{n+1}_{X}
\ar[r]_{ev}&X\ar@{=>}"a";"b"_{ad(F)}}$$ Here we use that $\S S^n_{\L
X}=S^{n+1}_X$ and $ad(ev(\S^{n-1}f))=ev(\S^n f)$.

The reader can check that the diagram
\begin{equation*}\tag{a}
\xymatrix{(\pi_{n,0}\L X)_{ab}\otimes(\pi_{n,0}\L
X)_{ab}\ar[r]^<(.3)\omega\ar@{=}[d]&
\pi_{n,1}\L X\ar[r]^\partial\ar[d]&\pi_{n,0}\L X\ar@{=}[d]\\
(\pi_{n+1,0}X)_{ab}\otimes(\pi_{n+1,0}X)_{ab}\ar[r]_<(.3)\omega&
\pi_{n+1,1}X\ar[r]_\partial&\pi_{n+1,0}X}
\end{equation*}
commutes, so it is a morphism of stable quadratic modules. Moreover,
the following diagram commutes
\begin{equation*}\tag{b}
\xymatrix{\pi_{n+1}\L X\ar@{^{(}->}[r]^\iota\ar[d]_{-\L}^\cong&
\pi_{n,1}\L X\ar[r]^\partial\ar[d]&\pi_{n,0}\L X\ar@{=}[d]\ar@{->>}[r]^ q &\pi_n\L X\ar[d]^\cong_\L\\
\pi_{n+2}X\ar@{^{(}->}[r]_\iota&
\pi_{n+1,1}X\ar[r]_\partial&\pi_{n+1,0}X\ar@{->>}[r]_ q &\pi_{n+1}X}
\end{equation*}
Here the exact rows are given by Proposition \ref{exacta} and the
arrows with $\cong$ are (up to sign) the usual isomorphisms of
homotopy groups, therefore the central vertical arrow in (a) is an
isomorphism by the five lemma.

For $n=2$ we have $\pi_{2,0}\L X=\grupo{\L^{3}X}_\ni=\pi_{3,0}X$ and
there is a homomorphism $\pi_{2,1}\L X\r\pi_{3,1}X$ defined as
above. This homomorphism makes commutative diagrams (a) and (b),
therefore it defines an isomorphism of reduced quadratic modules.

For $n=1$ there is an obvious epimorphism $\pi_{1,0}\L
X=\grupo{\L^2X}\twoheadrightarrow\grupo{\L^2X}_\ni=\pi_{2,0}X$. One
can also define a homomorphism $\pi_{1,1}\L X\r\pi_{2,1}X$ as above. It is easy to check that the
following square defines the desired crossed module morphism
\begin{equation*}\tag{c}
\xymatrix{
\pi_{1,1}\L X\ar[r]^\partial\ar[d]&\pi_{1,0}\L X\ar@{->>}[d]\\
\pi_{2,1}X\ar[r]_\partial&\pi_{2,0}X}
\end{equation*}
Moreover, the following diagram commutes
\begin{equation*}
\xymatrix{\pi_{2}\L X\ar@{^{(}->}[r]^\iota\ar[d]_{-\L}^\cong&
\pi_{1,1}\L X\ar[r]^\partial\ar[d]&\pi_{1,0}\L X\ar@{->>}[d]\ar@{->>}[r]^ q &\pi_1\L X\ar[d]^\cong_\L\\
\pi_{3}X\ar@{^{(}->}[r]_\iota&
\pi_{2,1}X\ar[r]_\partial&\pi_{2,0}X\ar@{->>}[r]_ q &\pi_{2}X}
\end{equation*}
This diagram is the analogue to (b) and shows that (c) is a weak
equivalence.

Now for $n=0$ we define the functor $\pi_{0,*}\L
X\r\phi_1\pi_{1,*}X$. On objects it is given by the inclusion
$Ob\pi_{0,*}\L X=\L X\subset\grupo{\L X}=Ob\phi_1\pi_{1,*}X$. Given
any object $f\in \L X$ in $\pi_{0,*}\L X$ we consider the inclusion
$\bar{f}\colon S^1\r S^1_X$ of the factor of the coproduct $S^1_X$
corresponding to $f$. Clearly the adjoint $ad(\bar{f}_{ev})\colon
S^0\r\L X$ is the inclusion of the point $f\in\L X$. If $g\in\L X$
is another object then a morphism $H\colon f\r g$ in $\pi_{0,*}\L X$
is just a track $H\colon ad(\bar{f}_{ev})\rr ad(\bar{g}_{ev})$ in
$\C{Top}^*$. The functor sends the morphism $H$ to the element in
$\pi_{1,0}X\ltimes\pi_{1,1}X$, which is $\pi_{1,0}X\times\pi_{1,1}X$
as a set, with $(\pi_1\bar{f})(1)$ in the left coordinate and
right coordinate given by the map
$$S^1\st{\mu}\To S^1\vee S^1\st{\nu\vee 1}\To S^1\vee S^1\st{(\bar{g},\bar{f})}\To S^1_X$$
and the track
$$\xymatrix{&&&S^1_X\ar@/^10pt/[rd]^{ev}&\\S^1\ar[r]_{\mu}\ar@/^20pt/[rrru]^0_{\;}="d"&
S^1\vee S^1\ar[r]_{\nu\vee 1}^{\;}="c"& S^1\vee
S^1\ar[r]_{(\bar{f},\bar{g})}^<(.6){\;}="a"\ar[ru]^{(\bar{g},\bar{g})}_<(.63){\;}="b"&
S^1_X\ar[r]_{ev}&
X\ar@{=>}"a";"b"_{(ad(H),0^\vc)}\ar@{=>}"c";"d"_N}$$ Here $N$ is a
nil-track and $ad(H)\colon \bar{f}_{ev}\rr \bar{g}_{ev}$ is the
adjoint of the track $H$. We leave it to the reader to check that
$\pi_{0,*}\L X\r\phi_1\pi_{1,*}X$ is a well-defined functor. One can
use again Proposition \ref{exacta} to check that this functor is an
equivalence.
\end{proof}

The functors $\phi_n$ in (\ref{loops1}) have left adjoints
\begin{equation}\label{losadjuntos}
\ad_n\colon\C{group}(n-1)\To\C{group}(n)
\end{equation}
given by $\ad_n=1_\C{Ab}$ for $n\geq 3$,
$$\ad_2\colon\C{Gr}\To\C{Ab},\text{ the abelianization};$$
$$\ad_1\colon\C{Set}^*\To\C{Gr},\text{ taking free group}.$$
These adjoints can be used to define for a pointed space $X$ the
natural suspension morphisms
$$\S\colon\ad_{n+1}\pi_nX\To\pi_{n+1}\S X$$
as the adjoint of
$$\pi_nX\st{\pi_nad(1)}\To\pi_{n}\L\S X\st{\L}\cong\phi_{n+1}\pi_{n+1}\S X.$$
Here we use the map $ad(1)\colon X\r\L\S X$ which is adjoint to the
identity in $\S X$ and the natural isomorphism $\L$. Now we generalize the situation for secondary
homotopy groups.

The functors $\phi_n$ in (\ref{loops2}) have left adjoints,
\begin{equation}\label{losadjuntos2}
\ad_n\colon\C{cross}(n-1)\To\C{cross}(n)
\end{equation}
given by $\ad_n=1_\C{squad}$ if $n\geq4$,
\begin{equation}\label{susps}
\begin{array}{c}
\ad_3\colon\C{rquad}\To\C{squad},\\
\ad_2\colon\C{cross}\To\C{rquad},\\
\ad_1\colon\C{grd}^*\To\C{cross}.
\end{array}
\end{equation}

\begin{lem}\label{haywe} 
The functors in (\ref{susps}) preserve $0$-free objects and weak equivalences
between them. 
\end{lem}

Here by convention we set all pointed groupoids to be $0$-free.

\begin{proof}[Proof of Lemma \ref{haywe}]
For $\ad_1$ the lemma follows from Lemma \ref{susex} below.

For $\ad_2$ and $\ad_3$ the lemma follows from the technical fact that the suspension  
 functors between
crossed and quadratic complexes described in \cite{ch4c} and \cite{cosqc}, 
which are extensions of $\ad_2$ and $\ad_3$, are compatible with the homotopy
relation in the category of totally free (i.e. cofibrant) crossed or quadratic complexes.
In addition we use that $0$-free crossed or quadratic modules are exactly the truncations of totally free
crossed or quadratic complexes.
\end{proof}

The functor $\ad_3$ is the stabilization in \cite{ch4c}
IV.C.3. It is defined as follows. Given a reduced quadratic module
$$(\omega,\partial)=(\otimes^2N_\abb\st{\omega}\To M\st{\partial}\To N)$$
the stabilized stable quadratic module 
$$\ad_3(\omega,\partial)=(\otimes^2N_\abb\st{\omega_\S}\To M_\S\st{\partial_\S}\To N)$$
is given by the group $M_\S$ obtained by quotienting out in $M$ the relations $$\omega(a\otimes b+b\otimes
a),\;\; a,b\in N_\abb,$$
and the homomorphisms $\omega_\S$ and $\partial_\S$ are induced by $\omega$ and $\partial$, respectively, in the
obvious way.

 The  functor $\ad_2$ in (\ref{susps}) is the
suspension functor in \cite{cosqc} 3.3. Given a crossed module $\partial\colon M\r N$ the reduced quadratic
module 
$$\ad_2\partial=(\otimes^2N_\abb\st{\omega}\To M^{\tilde{\S}}\st{\delta}\To N_\ni)$$
is given by the group
$M^{\tilde{\S}}$ which is a quotient of $M\times (\otimes^2N_\abb)$ by the relations
\begin{equation*}
(-m+m^n,0)=(0,\set{\partial(m)}\otimes\set{n})=(0,-\set{n}\otimes\set{\partial(m)}),
\end{equation*}
for any $m\in M$ and $n\in N$; and the homomorphisms $\delta$ and $\omega$ are defined by the following formulas,
$m\in M$, $n,n'\in N$,
\begin{eqnarray*}
\delta(m,\set{n}\otimes\set{n'})&=&\partial(m)+[n,n'],\\
\omega(\set{n}\otimes\set{n'})&=&(0,\set{n}\otimes\set{n'}).
\end{eqnarray*}

Finally we describe the functor $\ad_1$. Let $\C{G}$ be a groupoid with object pointed set $Ob\C{G}$ and morphism set $Mor\C{G}$. The crossed module $\ad_1\C{G}$ is the quotient of the free crossed module, see
\cite{ch4c}, generated by the function
$$Mor\C{G}\To\grupo{Ob\C{G}}$$
$$(h\colon U\r V)\mapsto -U+V$$
by the relations $u=g+f$ for $u,f,g\in Mor\C{G}$ with $u=fg$, the composition of $f$ and $g$ in $\C{G}$. One
readily checks that $\ad_1$ is the adjoint of $\phi_1$.

The functor $h_0$ commutes with $\ad_n$
\begin{equation}
h_0\ad_n=\ad_n h_0,\;\; n\geq 1.
\end{equation}
This follows from the definition of $\ad_n$ above for $n\geq 2$ and from Lemma \ref{susex} below in case $n=1$.
For $h_1$ the corresponding commutativity law is not true in general, compare Lemma \ref{susex} below.

\begin{thm}\label{sesusp}
There are natural morphisms in $\C{cross}(n+1)$
$$\Sigma\colon\ad_{n+1}\pi_{n,*}X\To\pi_{n+1,*}\Sigma X,\;\;n\geq0,$$
which induce the classical suspension homomorphism
$\S\colon\ad_{n+1}\pi_nX\r\pi_{n+1}\S X$ in $h_0$, and for $n\geq 3$ the homomorphism
$-\S\colon\ad_{n+2}\pi_{n+1}X\r\pi_{n+2}\S X$ in $h_1$. Moreover,
for $n\geq 3$ the morphism $\S$ is a weak equivalence provided $X$
is $m$-connected and $n\leq 2m-1$. It is also a weak equivalence for $n=2$ provided $X$ is simply connected, and
for
$n=1$ if  $X$ is connected. Furthermore, $\S$ is always a weak
equivalence for $n=0$.
\end{thm}

In the proof of this theorem we will use the following lemma.

\begin{lem}\label{susex}
For any pointed groupoid $\C{G}$ there are natural isomorphisms
\begin{enumerate}
\item $h_0\ad_1\C{G}=\grupo{\iso\C{G}}$,
\item
$h_1\ad_1\C{G}=\bigoplus\limits_{x\in\iso(\C{G})}(\aut_\C{G}(x))_\abb\otimes R$.
\end{enumerate}
Here $R$ is the group ring of $\grupo{\iso\C{G}}$.
\end{lem}

\begin{proof}
The crossed module $\ad_1\C{G}$ defined above is the
truncation $$N_1FB\C{G}/d(N_2FB\C{G})\r N_0FB\C{G}$$ of the Moore
complex $N_*FB\C{G}$ of the Milnor construction $FB\C{G}$ on the
classifying space $B\C{G}$ of the pointed groupoid $\C{G}$, see \cite{cdhg}
and \cite{gj} I.1.4 and V.6. To see this we have on the $0$-level 
$$N_0FB\C{G}=\grupo{Ob\C{G}},$$
and on the $1$-level the set $Mor\C{G}$ in $\ad_1\C{G}$ is mapped to $N_1FB\C{G}/d(N_2FB\C{G})$
by sending $h\colon U\r V$ to the coset modulo $d(N_2FB\C{G})$ of the element
$$-1_U+h\in N_1FB\C{G}\subset F_1B\C{G}=\grupo{Mor\C{G}}.$$
This can be checked by computing $N_1FB\C{G}/d(N_2FB\C{G})$ in terms of generators and relations. In order to
carry out this computation one uses the Reidemeister-Schreier method, see \cite{mks}, 
which simplifies in this particular
case since the simplicial identities hold in $FB\C{G}$ and the boundaries and degeneracies in this
simplicial group are homomorphisms between free groups on pointed sets induced by maps between the generating
pointed sets.

By the previous observation the kernel and cokernel of
$\ad_1\C{G}$ are the $\pi_1$ and $\pi_2$ of the suspension $\S|B\C{G}|$
of the geometric realization $|B\C{G}|$ of the classifying space of
$\C{G}$, so the lemma follows from elementary facts from homotopy
theory.

We also remark that given $x\in\iso(\C{G})$ and a representative
$\tilde{x}\in Ob\C{G}$ of $x$ the group
$(\aut_\C{G}(\tilde{x}))_\abb$ does not depend on the choice of
$\tilde{x}$, up to natural isomorphism, therefore we can denote it
by $(\aut_\C{G}(x))_\abb$.
\end{proof}

\begin{proof}[Proof of Theorem \ref{sesusp}]
Consider the morphism
$$\pi_{n,*}X\st{\pi_{n,*}ad(1)}\To\pi_{n,*}\L\S X\st{\L}\To\phi_{n+1}\pi_{n+1,*}\S X,$$ where $ad(1)\colon X\r\L\S X$ is the adjoint of
the identity in $\S X$ and $\L$ is given by Proposition \ref{adeq}.
The morphism in the statement is the adjoint
of this one. 
For $n\geq 3$ the range where this morphism is a weak equivalence
follows from Proposition \ref{exacta} and the classical suspension
theorem for ordinary homotopy groups.

For $n=1$ the theorem follows from Proposition \ref{concel} below
and \cite{cosqc} 4.8. For $n=2$ we use Proposition \ref{concel} and
\cite{ch4c} IV.C. For this we use that we are dealing with $0$-free
objects and that $\ad_n$ preserves weak equivalences between them, see Lemma \ref{haywe}.

If $n=0$ we have $\iso(\pi_{0,*}X)=\pi_0X$ and for any $x\in
Ob\pi_{0,*}X$, $\aut_{\pi_{0,*}X}(x)=\pi_1(X,x)$. By using
elementary homotopy theory one can check that $$\pi_1\S
X\cong\grupo{\pi_0X}$$ and
$$\pi_2\S
X\cong\bigoplus\limits_{x\in\pi_0X}(\pi_1(X,x))_\abb\otimes\Z\grupo{\pi_0X}.$$
Now it is enough to notice that isomorphisms in Lemma \ref{susex}
are compatible with the two isomorphisms above and Proposition
\ref{exacta} (in this last case up to sign $-1$ in kernel).
\end{proof}

\section{Secondary homotopy groups as $2$-functors}\label{functor}

In Section \ref{2hg} we have defined the secondary homotopy group
functors $(n\geq 0)$ $$\pi_{n,*}\colon\C{Top}^*\To\C{cross}(n).$$ As
we recalled in Section \ref{tbm} the category of pointed spaces is a
groupoid-enriched category, therefore it is reasonable to wonder whether
$\pi_{n,*}$ is a $2$-functor. This is known to be true if $n=0$.
In this section we prove that it is actually true for any $n\geq 0$.

We recall from \cite{ch4c} the definition of homotopies or $2$-morphisms in the categories of crossed modules and reduced quadratic modules.

\begin{defn}\label{tcq}

Suppose that we have two crossed modules $\partial\colon M\r N$, ${\partial'}\colon {M'}\r{N'}$ and two morphisms
$f,g\colon\partial\r{\partial'}$ given by
$$\xymatrix{M\ar[d]_\partial\ar[r]^{f_1,g_1}&M'\ar[d]^{\partial'}\\N\ar[r]_{f_0,g_0}&N'}$$
A \emph{$2$-morphism} $\alpha\colon f\rr g$ is a 
function $\alpha\colon N\r M'$ such that for any $x,y\in N$ and any $m\in M$ the following equalities hold:
\begin{enumerate}
\item $\alpha(x+y)=\alpha(x)^{f_0(y)}+\alpha(y)$,
\item $g_0(x)=f_0(x)+\partial'(\alpha(x))$,
\item $g_1(m)=f_1(m)+\alpha(\partial(m))$.
\end{enumerate}

If we now have two reduced quadratic modules
$$\otimes^2N_\abb\st{\omega}\To M\st{\partial}\To N,$$
$$\otimes^2(N')_\abb\st{\omega'}\To M'\st{\partial}\To N',$$
and two morphisms $f,g\colon(\omega,\partial)\r(\omega',\partial')$
$$\xymatrix{\otimes^2 N_\abb\ar[r]^\omega\ar[d]^{\otimes^2 (f_0)_\abb}_{\otimes^2 (g_0)_\abb}&M\ar[r]^\partial\ar[d]^{f_1}_{g_1}&N\ar[d]^{f_0}_{g_0}\\
\otimes^2(N')_\abb\ar[r]_{\omega'}&M'\ar[r]_{\partial'}&N'}$$ a
$2$-morphism $\alpha\colon f\rr g$ is just a $2$-morphism $\alpha\colon\phi_2 f\rr
\phi_2 g$ in the category of crossed modules. Here we use the
forgetful functor $\phi_2$ in (\ref{loops}). More precisely,
$\alpha\colon N\r M'$ is a function such that for any $x,y\in N$ and
$m\in M$ the following equations hold:
\begin{enumerate}
\item $\alpha(x+y)=\alpha(x)+\alpha(y)+\omega'(\set{-f_0(x)+g_0(x)}\otimes\set{f_0(y)})$,
\item $g_0(x)=f_0(x)+\partial'(\alpha(x))$,
\item $g_1(m)=f_1(m)+\alpha(\partial(m))$.
\end{enumerate}

The $2$-morphisms for stable quadratic module morphisms are the same as $2$-morphisms
for the corresponding reduced quadratic module morphisms. In
particular the forgetful functors $\phi_n$ in (\ref{loops2}) become
automatically $2$-functors which are full and faithful at the level
of $2$-morphisms for $n\geq 2$.

The $2$-morphisms in the category $\C{grd}^*$ of pointed groupoids are just
natural transformations between functors. For $n=0$ the functor
$\phi_1$ in (\ref{loops}) is also a $2$-functor. More precisely,
if $\partial\colon M\r N$ and $\partial'\colon M'\r N'$ are crossed
modules and $\alpha\colon N\r M'$ is a $2$-morphism $\alpha\colon f\rr g$
between two morphisms $f,g\colon\partial\r\partial'$ then the
natural transformation $\phi_1\alpha\colon\phi_1f\rr\phi_1g$ between
the pointed groupoid morphisms
$\phi_1f,\phi_1g\colon\phi_1\partial\r\phi_1\partial'$ is given by
the morphisms $(f_0(n),\alpha(n))\colon f_0(n)\r g_0(n)$ in
$\phi_1\partial'$ which are natural in $n\in N$.
\end{defn}

\begin{prop}\label{estc}
The category $\C{cross}(n)$ with $2$-morphisms as in Definition \ref{tcq}
is a groupoid-enriched category. Moreover, for all $n\geq 1$ the functor $\ad_n$ in (\ref{susps}) is a $2$-functor  which is adjoint to
$\phi_n$ in (\ref{loops}) as a groupoid-enriched functor.
\end{prop}

\begin{proof}
For $n=0$ it is well-known that $\C{cross}(n)$ is a groupoid-enriched category. We only need to carry out the
proof of the first part of the statement for crossed modules since the groupoid-enriched structure in $\C{rquad}$
and $\C{squad}$ is pulled back through the forgetful functors in
(\ref{loops}).

In this proof  $\partial_i\colon M_i\r N_i=\grupo{E_i}$ will denote a crossed module for  $i=0,1,2,3$.

Let $f,g,h\colon \partial_1\r\partial_2$ be crossed module morphisms and let $\alpha\colon f\rr g$, $\beta\colon
g\rr h$ be vertically composable $2$-morphisms. The vertical composition is defined by
$(\beta\vc\alpha)(x)=\alpha(x)+\beta(x)$ for any $x\in M_1$. The inverse $2$-morphism $\alpha^\vi\colon g\rr f$ is
defined by $\alpha^\vi(x)=-\alpha(x)$ and the trivial $2$-morphism $0^\vc_f\colon f\rr f$ is $0^\vc_f(x)=0$.

Suppose that we have a diagram
\begin{equation*}
\xymatrix@C=35pt{\partial_0\ar[r]^{k}&\partial_1\ar@/_15pt/[r]_{g}^{\;}="a"
\ar@/^15pt/[r]^{f}_{\;}="b"&\partial_2
\ar[r]^{h}&\partial_3\ar@{=>}"b";"a"^\alpha}
\end{equation*}
Then the two possible horizontal compositions $\alpha k\colon fk\rr gk$, $h\alpha\colon hf\rr gf$ are defined as $\alpha k=\alpha k_0\colon N_0\r M_3$ and
$h\alpha=h_1\alpha\colon N_1\r M_4$.

Suppose now that we have a diagram
\begin{equation*}
\xymatrix@C=35pt{\partial_0\ar@/_10pt/[r]_g^{\;}="a"\ar@/^10pt/[r]^{f}_{\;}="b"&\partial_1
\ar@/_10pt/[r]_{g'}^{\;}="c"\ar@/^10pt/[r]^{f'}_{\;}="d"&\partial_2\ar@{=>}"b";"a"^\alpha\ar@{=>}"d";"c"^{\alpha'}}
\end{equation*}
Let us check the equality
\begin{equation*}\tag{a}
(g'\alpha)\vc(\alpha'f)=(\alpha'g)\vc(f'\alpha).
\end{equation*}
Given $x\in M_0$ by using the equations defining crossed modules and $2$-morphisms we get
\begin{eqnarray*}
((g'\alpha)\vc(\alpha'f))(x)&=&\alpha'f_0(x)+g'_1\alpha(x)\\
&=&\alpha'f_0(x)+f'_1\alpha(x)+\alpha'\partial_1\alpha(x)\\
&=&\alpha'f_0(x)+f'_1\alpha(x)+\alpha'(-f_0(x)+g_0(x))\\
&=&\alpha'f_0(x)+f'_1\alpha(x)+\alpha'(-f_0(x))^{f'_0g_0(x)}+\alpha'(g_0(x))\\
&=&\alpha'f_0(x)+f'_1\alpha(x)-\alpha'(f_0(x))^{f'_0(-f_0(x)+g_0(x))}+\alpha'(g_0(x))\\
&=&\alpha'f_0(x)+f'_1\alpha(x)-\alpha'(f_0(x))^{f'_0\partial_1\alpha(x)}+\alpha'(g_0(x))\\
&=&\alpha'f_0(x)+f'_1\alpha(x)-\alpha'(f_0(x))^{\partial_2f'_1\alpha(x)}+\alpha'(g_0(x))\\
&=&\alpha'f_0(x)-\alpha'(f_0(x))+f'_1\alpha(x)+\alpha'(g_0(x))\\
&=&f'_1\alpha(x)+\alpha'(g_0(x))\\
&=&((\alpha'g)\vc(f'\alpha))(x).
\end{eqnarray*}
Hence (a) holds and $\C{cross}(1)$ is indeed a groupoid-enriched category.

Now we define the $2$-functors $\ad_n$ at the level of $2$-morphisms. For $n\geq 4$ they are identity $2$-morphism
functors. For $n=3$, given a $2$-morphism $\alpha\colon N\r M'$ between two reduced quadratic module morphisms from
$(\omega,\partial)$ to $(\omega',\partial')$ the $2$-morphism $\ad_3\alpha\colon N\r M'_\S$ is the composition of
$\alpha$ with the natural projection $M'\twoheadrightarrow M'_\S$. For $n=2$, if $\alpha\colon N\r M'$ is a $2$-morphism
between two crossed module morphisms from $\partial$ to $\partial'$ then $\ad_2\alpha\colon N_\ni\r
(M')^{\tilde{\S}}$ is defined as $(\ad_2\alpha)(n)=(\alpha(n),0)$ for $n\in N_\ni$. Finally for $n=1$, if
$\alpha\colon f\rr g$ is a natural transformation between two pointed groupoid morphisms $f,g\colon\C{G}\r\C{G}'$
which is given by a collection of morphisms $\alpha(X)\colon f(X)\r g(X)$ in $\C{G}'$ for $X\in Ob\C{G}$
then $\ad_1\alpha\colon \grupo{Ob\C{G}}\r(\ad_1\C{G}')_1$ is the unique $2$-morphism between crossed module morphisms
satisfying $(\ad_1\alpha)(X)=\alpha(X)$.
\end{proof}

\begin{thm}\label{aretrack}
The secondary homotopy groups are $2$-functors $(n\geq 0)$
$$\pi_{n,*}\colon\C{Top}^*\To\C{cross}(n).$$
\end{thm}

In the proof of Theorem \ref{aretrack} we will use the following
general construction.

\begin{defn}\label{rcon}
Let $X$ be a pointed space. Given two maps $f,g\colon S^1\r
\vee_{\L^nX}S^1$ and a track $H\colon f_{ev}\rr g_{ev}$ we define
$r(H)\in\pi_{n,1}X$ as follows. Let $\varepsilon\colon S^1\r S^1\vee
S^1$ be a map with $\pi_1\varepsilon\colon\Z\r\grupo{a,b}$ satisfying
$(\pi_1\varepsilon)(1)=-a+b$, or just
$(\pi_1\varepsilon)_\ni(1)=-a+b$ if $n\geq 2$, then $r(H)$ is
represented by the map
$$S^1\st{\varepsilon}\To S^1\vee S^1\st{(f,g)}\To\vee_{\L^nX}S^1$$
and the track
$$\xymatrix{&S^n\ar@/^15pt/[rrd]^{g_{ev}}_<(.2){\;}="d"&&\\
S^n\ar[r]_{\S^{n-1}\varepsilon}^{\;}="a"\ar@/^10pt/[ru]^0_<(.2){\;}="b"&
S^n\vee S^n\ar[u]^{(1,1)}\ar[r]_{\S^{n-1}(f,g)}^{\;}="c"&S^n_X\ar[r]_{ev}&X\ar@{=>}"a";"b"_N\ar@{=>}"c";"d"_{(H,0^\vc)}}$$
\end{defn}

The $r$-construction may be regarded as a machine to generate $2$-morphisms in $\C{cross}(n)$ between secondary
homotopy groups. Some of the axioms of a $2$-morphism in $\C{cross}(n)$ are checked in the following lemma for
the $r$-construction.

\begin{lem}\label{general}
Let $X$ be a pointed space. Given $f,g,h\colon S^1\r \vee_{\L^nX}S^1$, $H\colon f_{ev}\rr g_{ev}$ and $K\colon g_{ev}\rr h_{ev}$
the following formulas hold in $\pi{n,1}X$,
\begin{enumerate}
\item $\partial r(H)= -(\pi_1f)(1)+(\pi_1g)(1)$ if $n=1$,
\item $\partial r(H)= -(\pi_1f)_\ni(1)+(\pi_1g)_\ni(1)$ if $n\geq 2$,
\item $r(K\vc H)= r(H)+r(K)$.
\end{enumerate}
\end{lem}

\begin{proof}
Equations (1) and (2) are clear. The element $r(H)+r(K)$ is represented by
the following diagram
\begin{equation*}\tag{a}
\xymatrix{&&S^n\vee S^n\ar@/^13pt/[rrrd]^{(g_{ev},h_{ev})}_<(.2){\;}="d"&&&
\\S^n\ar[r]_-{\S^{n-1}\mu}&S^n\vee
S^n\ar@/^15pt/[ru]^0_<(.2){\;}="b"\ar[r]_-{\S^{n-1}(\varepsilon\vee\varepsilon)}^{\;}="a"&
S^n\vee S^n\vee S^n\vee S^n\ar[u]^{(i_1,i_1,i_2,i_2)}\ar[rr]_-{\S^{n-1}(f,g,g,h)}^{\;}="c"&&S^n_X\ar[r]_-{ev}&X
\ar@{=>}"a";"b"_N\ar@{=>}"c";"d"_{(H,0^\vc,K,0^\vc)}}
\end{equation*}
By Remark \ref{*} we can insert the track $(K,0^\vc)$ in (a) so that the pasting of (a) and (b) below remains the same
\begin{equation*}\tag{b}
\xymatrix{&&S^n\ar@/^25pt/[rrrdd]^{h_{ev}}_<(.2){\;}="f"&&&\\
&&S^n\vee S^n\ar@/^15pt/[rrrd]^{(g_{ev},h_{ev})}_<(.2){\;}="d"^<(.2){\;}="e"\ar[u]^{(1,1)}&&&
\\S^n\ar[r]_-{\S^{n-1}\mu}_<(.7){\;}="g"&S^n\vee
S^n\ar@/^15pt/[ru]^0_<(.2){\;}="b"\ar[r]_-{\S^{n-1}(\varepsilon\vee\varepsilon)}^{\;}="a"&
S^n\vee S^n\vee S^n\vee S^n\ar[u]^{(i_1,i_1,i_2,i_2)}\ar[rr]_-{\S^{n-1}(f,g,g,h)}^{\;}="c"&&S^n_X\ar[r]_-{ev}&X
\ar@{=>}"a";"b"_N\ar@{=>}"c";"d"_{(H,0^\vc,K,0^\vc)}
\ar@{=>}"e";"f"_{(K,0^\vc)}}
\end{equation*}
Pasting some tracks in (b) we obtain
\begin{equation*}\tag{c}
\xymatrix@R=50pt{&&S^n\ar@/^20pt/[rrrd]^{h_{ev}}_<(.2){\;}="d"&&&
\\S^n\ar[r]_-{\S^{n-1}\mu}&S^n\vee
S^n\ar@/^20pt/[ru]^0_<(.2){\;}="b"\ar[r]_-{\S^{n-1}(\varepsilon\vee\varepsilon)}^{\;}="a"&
S^n\vee S^n\vee S^n\vee S^n\ar[u]^{(1,1,1,1)}\ar[rr]_-{\S^{n-1}(f,g,g,h)}^{\;}="c"&&S^n_X\ar[r]_-{ev}&X
\ar@{=>}"a";"b"_N\ar@{=>}"c";"d"_{(K\vc H,K,K,0^\vc)}}
\end{equation*}
One can factor (c) as
\begin{equation*}\tag{d}
\xymatrix@R=50pt{&&S^n\ar@/^20pt/[rrrd]^{h_{ev}}_{\;}="d"&&&
\\S^n\ar[r]_-{\S^{n-1}\mu}&S^n\vee
S^n\ar@/^20pt/[ru]^0_<(.2){\;}="b"\ar[r]_-{\S^{n-1}(\varepsilon\vee\varepsilon)}^{\;}="a"&
S^n\vee S^n\vee S^n\vee S^n\ar[u]^{(1,1,1,1)}\ar[r]_-{(i_1,i_2,i_2,i_3)}&
S^n\vee S^n\vee S^n\ar[ul]^{(1,1,1)}\ar[r]_-{\S^{n-1}(f,g,h)}^{\;}="c"&S^n_X\ar[r]_-{ev}&X
\ar@{=>}"a";"b"_N\ar@{=>}"c";"d"|{(K\vc H,K,0^\vc)}}
\end{equation*}
Diagram (d) represents the same element in $\pi_{n,1}X$ as
\begin{equation*}\tag{e}
\xymatrix@R=50pt{&&S^n\ar@/^20pt/[rrrd]^{h_{ev}}_{\;}="d"&&&
\\S^n\ar[r]^{\S^{n-1}\mu}_<(.7){\;}="f"\ar@/_20pt/[rrd]_-{\S^{n-1}\varepsilon}^<(.2){\;}="e"&S^n\vee
S^n\ar@/^20pt/[ru]^0_<(.2){\;}="b"\ar[r]_-{\S^{n-1}(\varepsilon\vee\varepsilon)}^{\;}="a"&
S^n\vee S^n\vee S^n\vee S^n\ar[u]^{(1,1,1,1)}\ar[r]_-{(i_1,i_2,i_2,i_3)}&
S^n\vee S^n\vee S^n\ar[ul]^{(1,1,1)}\ar[r]_-{\S^{n-1}(f,g,h)}^{\;}="c"&S^n_X\ar[r]_-{ev}&X
\ar@{=>}"a";"b"_N\ar@{=>}"c";"d"|{(K\vc H,K,0^\vc)}\ar@{=>}"e";"f"_N\\
&&S^n\vee S^n\ar@/_20pt/[ru]_-{(i_1,i_3)}}
\end{equation*}
where the two $N$ denote different nil-tracks. The pasting of (e) coincides with the pasting of
\begin{equation*}\tag{f}
\xymatrix{&S^n\ar@/^15pt/[rrd]^{h_{ev}}_<(.2){\;}="d"&&\\
S^n\ar[r]_-{\S^{n-1}\varepsilon}^{\;}="a"\ar@/^10pt/[ru]^0_<(.2){\;}="b"&
S^n\vee
S^n\ar[u]^{(1,1)}\ar[r]_-{\S^{n-1}(f,h)}^{\;}="c"&S^n_X\ar[r]_-{ev}&X\ar@{=>}"a";"b"_N\ar@{=>}"c";"d"_{(K\vc
H,0^\vc)}}
\end{equation*}
and (f) represents $r(K\vc H)$, hence (3) follows.
\end{proof}

In the next lemma we check the derivation property for the $r$-construction.

\begin{lem}\label{general2}
Let $X, Y$ be pointed spaces. Given $x,y\colon S^1\r
\vee_{\L^nX}S^1$,
$f,g\colon\vee_{\L^nX}S^1\r\vee_{\L^nY}S^1$, and $H\colon
ev(\S^{n-1}f)\rr ev(\S^{n-1}g),$ we have the following equalities
in $\pi_{n,1}Y$:
\begin{eqnarray*}
r(H(\S^{n-1}(x,y))(\S^{n-1}\mu))&=&r(H(\S^{n-1}x))^{\pi_1(fy)(1)}+r(H(\S^{n-1}y))\text{
if }n=1,\\
\text{and if }n\geq
2&=&r(H(\S^{n-1}x))+r(H(\S^{n-1}y))\\&&+\omega(\set{-\pi_1(fx)(1)+\pi_1(gx)(1)},\set{\pi_1(fy)(1)}).
\end{eqnarray*}
\end{lem}

\begin{proof}
Suppose that $n=1$. Let $\varpi\colon S^1\r S^1\vee \st{5}\cdots\vee
S^1$ be a map with
$$\pi_1\varpi\colon\Z\r\grupo{a_1,a_2,a_3,a_4,a_5}\colon 1\mapsto
-a_1-a_2+a_3+a_1-a_4+a_5.$$ The element
$r(H(\S^{n-1}x))^{\pi_1(fy)(1)}+r(H(\S^{n-1}y))$ is represented by
the diagram
\begin{equation*}\tag{a}
\xymatrix{&&S^n_X\vee S^n_X\ar[rr]^{\S^{n-1}(f,g)}_{\;}="d"&&S^n_Y\ar@/^10pt/[rd]^{ev}&
\\S^n\ar@/^13pt/[rru]^0_<(.6){\;}="b"\ar[r]_<(.17){\S^{n-1}\varpi}&\vee_5S^n\ar[r]_<(.2){z}^{\;}="a"&S^n_X\vee S^n_X\vee
S^n_X\ar[u]^{(i_1,i_2,i_2)}\ar[rr]_<(.4){\S^{n-1}(f,f,g)}^{\;}="c"&&S^n_Y\ar[r]_{ev}&Y
\ar@{=>}"a";"b"^N\ar@{=>}"c";"d"_{(0^\vc,H,0^\vc)}}
\end{equation*} 
where $z=\S^{n-1}(i_1y,i_2x,i_3x,i_2y,i_3y)$ and $N$ is a nil-track. By Remark \ref{*} we can add the track
$(H,0^\vc)$ to (a) so that the pasting of (a) and (b) below remains the same
\begin{equation*}\tag{b}
\xymatrix{&&S^n_X\ar[rr]^{\S^{n-1}g}_{\;}="f"&&S^n_Y\ar@/^10pt/[rdd]^{ev}&\\
&&S^n_X\vee
S^n_X\ar[u]^{(1,1)}\ar[rr]|{\S^{n-1}(f,g)}_{\;}="d"^{\;}="e"&&S^n_Y\ar@/^10pt/[rd]_{ev}&
\\S^n\ar@/^13pt/[rru]^0_<(.6){\;}="b"\ar[r]_<(.17){\S^{n-1}\varpi}&\vee_5 S^n\ar[r]_<(.2){z}^{\;}="a"&S^n_X\vee
S^n_X\vee
S^n_X\ar[u]^{(i_1,i_2,i_2)}\ar[rr]_<(.4){\S^{n-1}(f,f,g)}^{\;}="c"&&S^n_Y\ar[r]_{ev}&Y
\ar@{=>}"a";"b"^N\ar@{=>}"c";"d"_{(0^\vc,H,0^\vc)}\ar@{=>}"e";"f"_{(H,0^\vc)}}
\end{equation*}
Pasting some tracks in (b) we obtain
\begin{equation*}\tag{c}
\xymatrix{&&S^n_X\ar[rr]^{\S^{n-1}g}_{\;}="d"&&S^n_Y\ar@/^10pt/[rd]^{ev}&
\\S^n\ar@/^13pt/[rru]^0_<(.6){\;}="b"\ar[r]_<(.17){\S^{n-1}\varpi}&\vee_5 S^n\ar[r]_<(.2){z}^{\;}="a"&S^n_X\vee S^n_X\vee
S^n_X\ar[u]^{(1,1,1)}\ar[rr]_<(.4){\S^{n-1}(f,f,g)}^{\;}="c"&&S^n_Y\ar[r]_{ev}&Y
\ar@{=>}"a";"b"^N\ar@{=>}"c";"d"_{(H,H,0^\vc)}}
\end{equation*}
The pasting of (c) is the same as the pasting of 
\begin{equation*}\tag{d}
\xymatrix{&&S^n_X\ar[rr]^{\S^{n-1}g}_{\;}="d"&&S^n_Y\ar@/^10pt/[rd]^{ev}&
\\S^n\ar@/^13pt/[rru]^0_<(.6){\;}="b"\ar[r]_<(.17){\S^{n-1}\varpi}&\vee_5 S^n
\ar[r]_<(.2){z}^{\;}="a"& S^n_X\vee
S^n_X\ar[u]^{(1,1)}\ar[rr]_<(.4){\S^{n-1}(f,g)}^{\;}="c"&&S^n_Y\ar[r]_{ev}&Y
\ar@{=>}"a";"b"^N\ar@{=>}"c";"d"_{(H,0^\vc)}}
\end{equation*}
Diagram (d)
represents the same element in $\pi_{1,1}Y$ as
\begin{equation*}\tag{e}
\xymatrix{&&S^n_X\ar[rr]^{\S^{n-1}g}_{\;}="d"&&S^n_Y\ar@/^10pt/[rd]^{ev}&
\\S^n\ar[d]_{\S^{n-1}\varepsilon}\ar@/^13pt/[rru]^0_<(.6){\;}="b"\ar[r]|<(.3){\S^{n-1}\varpi}_{\;}="f"&
\vee_5 S^n\ar[r]_<(.2){z}^{\;}="a"& S^n_X\vee
S^n_X\ar[u]^{(1,1)}\ar[rr]_<(.4){\S^{n-1}(f,g)}^{\;}="c"&&S^n_Y\ar[r]_{ev}&Y
\ar@{=>}"a";"b"^N\ar@{=>}"c";"d"_{(H,0^\vc)}\\S^n\vee
S^n\ar[r]_<(.2){\S^{n-1}(\mu\vee\mu)}^{\;}="e"&\vee_4
S^n\ar@/_13pt/[ru]_<(.4){\;\;\;\;\;\S^{n-1}(x,y)\vee(x,y)}\ar@{=>}"e";"f"_N}
\end{equation*}
Here the two $N$ denote different nil-tracks.
Diagram (e) represents $r(H(\S^{n-1}(x,y))(\S^{n-1}\mu))$, hence the case $n=1$ follows.

Suppose now that $n\geq 2$. By using Theorem \ref{propi}  and the
claim (*) in the proof of Proposition \ref{redquades} the element
$$r(H(\S^{n-1}x))+r(H(\S^{n-1}y))+\omega(\set{-\pi_1(fx)(1)+\pi_1(gx)(1)},\set{\pi_1(fy)(1)})$$
is represented by the diagram
\begin{equation*}\tag{f}
\xymatrix@C=28pt{&&&&S^n_X\ar[r]^{\S^{n-1}g}_{\;}="d"&S^n_Y\ar@/^10pt/[rd]^{ev}&
\\S^n\ar[d]^{\S^{n-1}\varepsilon}\ar@/^13pt/[rrrru]^0_<(.3){\;}="b"\ar[r]|<(.3){\S^{n-1}\mu}_<(.565){\;}="f"&
S^n\vee
S^n\ar[r]_<(.3){\S^{n-1}(\varepsilon\vee\varepsilon)}^{\;}="a"&\vee_4
S^n\ar[rr]_<(.35){\S^{n-1}(x\vee x,y\vee y)}&& S^n_X\vee
S^n_X\ar[u]^{(1,1)}\ar[r]_<(.4){\S^{n-1}(f,g)}^{\;}="c"&S^n_Y\ar[r]_{ev}&Y
\ar@{=>}"a";"b"_N\ar@{=>}"c";"d"_{(H,0^\vc)}\\
S^n\vee S^n\ar[rr]^<(.2){\;}="e"_{\S^{n-1}(\mu\vee\mu)}&&\vee_4
S^n\ar@/_15pt/[rru]_{\;\;\;\;\;\;\S^{n-1}(x,y)\vee(x,y)}\ar@{=>}"e";"f"_Q}
\end{equation*}
where $Q$ is the unique nil-track with Hopf invariant
$$-(-i_1(\pi_1(x))_\abb(1)+i_2(\pi_1(x))_\abb(1))\otimes(i_1(\pi_1(y))_\abb(1))\in\otimes^2_n(\Z[\L^nX]\oplus\Z[\L^nX]),$$
for $i_1,i_2\colon\Z[\L^nX]\r\Z[\L^nX]\oplus\Z[\L^nX]$ the inclusion
of the factors of the direct sum. By Theorem \ref{propi} the track
$(1,1)Q$ is a nil-track, hence diagram (f) also represents the
element $r(H(\S^{n-1}(x,y))(\S^{n-1}\mu))$, and we are done.
\end{proof}

Now we are ready to prove the main result of this section.

\begin{proof}[Proof of Theorem \ref{aretrack}]
This is known for $n=0$. Suppose now that $n\geq 1$. Let $f,g\colon
X\r Y$ be two maps and $f_*,g_*\colon\pi_{n,*}X\r\pi_{n,*}Y$ the
induced morphisms in $\C{cross}(n)$. Moreover, let $H\colon f\rr g$
be a track. We define a $2$-morphism $\pi_{n,*}H=H_*\colon f_*\rr g_*$ in
the following way. The function $H_*\colon\pi_{n,0}X\r\pi_{n,1}Y$
sends an element $x\in\pi_{n,0}X$ represented by a map
$\tilde{x}\colon S^1\r\vee_{\L^nX}S^1$ with $(\pi_1\tilde{x})(1)=x$
 if $n=1$ and
$(\pi_1\tilde{x})_\ni(1)=x$ if $n\geq 2$ to the element in
$\pi_{n,1}Y$ represented by the map (a) and the track (b)
below.

Let $\varepsilon\colon S^1\r S^1\vee S^1$ be a map with
$(\pi_1\varepsilon)(1)=-a+b\in\grupo{a,b}$. The map (a) is defined by
\begin{equation*}\tag{a}
\xymatrix{S^1\ar[r]^<(.3){\varepsilon}&S^1\vee
S^1\ar[r]^<(.15){\tilde{x}\vee\tilde{x}}&(\vee_{\L^nX}
S^1)\vee(\vee_{\L^nX} S^1) \ar[rr]^<(.35){\S\L^n(f
,g)}&&\vee_{\L^nY} S^1}.
\end{equation*} 

The track (b) is given by
\begin{equation*}\tag{b}
\xymatrix{&S^n\ar[rr]^{\S^{n-1}\tilde{x}}&&S^n_X\ar[rr]^{ev}&&X\ar@/^20pt/[rdd]^g_{\;}="b"&\\&&&&&X\vee
X\ar[rd]_{(f,g)}^<(.1){\;}="a"\ar[u]^{(1,1)}&
\\S^n\ar@/^20pt/[ruu]^0_{\;}="d"
\ar[r]_<(.3){\S^{n-1}\varepsilon}&S^n\vee
S^n\ar[rr]_{\S^{n-1}(\tilde{x}\vee\tilde{x})}\ar[uu]_{(1,1)}^{\;}="c"&&
S^n_X\vee S^n_X \ar[uu]_{(1,1)}\ar[rru]^{ev\vee
ev}\ar[rr]_<(.35){\S^n\L^n(f ,g )}&&
S^n_Y\ar[r]_{ev}&Y\ar@{=>}"a";"b"_<(.2){(H,0^\vc)}\ar@{=>}"c";"d"^N}
\end{equation*}
Here $N$ is a nil-track. With the terminology introduced in Definition \ref{rcon} we have $H_*(x)=r(H\, ev (\S^{n-1}\tilde{x}))$ for the
track $H\, ev \, (\S^{n-1}\tilde{x})\colon ev(\S^n\L^n f)(\S^{n-1}
\tilde{x})= f\, ev(\S^{n-1}\tilde{x})\rr g\, ev(\S^{n-1}\tilde{x}=ev(\S^n\L^ng)(\S^{n-1}\tilde{x}$.

The proof of equation (1) in the definition of $2$-morphisms in
$\C{cross}(n)$ follows from Lemma \ref{general2}. Equation (2)
follows from Lemma \ref{general} (1) or (2). Equation (3) follows
from the fact that given $[k,K]\in\pi_{n,1}X$ the pasting of the following
tracks coincide. 

The track (c) below represents $H_*\partial[k,K]$
\begin{equation*}\tag{c}
\xymatrix{&S^n\ar[rr]^{\S^{n-1}k}&&S^n_X\ar[rr]^{ev}&&X\ar@/^20pt/[rdd]^g_{\;}="b"&\\&&&&&X\vee X\ar[rd]_{(f,g)}^<(.1){\;}="a"\ar[u]^{(1,1)}&
\\S^n\ar@/^20pt/[ruu]^0_{\;}="d"
\ar[r]_<(.3){\S^{n-1}\varepsilon}&S^n\vee S^n\ar[rr]_{\S^{n-1}(k\vee
k)}\ar[uu]_{(1,1)}^{\;}="c"&& S^n_X\vee S^n_X
\ar[uu]_{(1,1)}\ar[rru]^{ev\vee ev}\ar[rr]_<(.35){\S^n\L^n(f
,g)}&&
S^n_Y\ar[r]_{ev}&Y\ar@{=>}"a";"b"_<(.2){(H,0^\vc)}\ar@{=>}"c";"d"^N}
\end{equation*}
Using Remark \ref{*} we can introduce the track $K$ on the top of (c) without altering the result of the pasting
\begin{equation*}\tag{d}
\xymatrix{&S^n\ar@/^30pt/[rrrr]^0_<(.7){\;}="f"\ar[rr]^{\S^{n-1}k}&&S^n_X\ar[rr]_{ev}^{\;}="e"
&&X\ar@/^20pt/[rdd]^g_{\;}="b"&\\&&&&&X\vee
X\ar[rd]_{(f,g)}^<(.1){\;}="a"\ar[u]^{(1,1)}&
\\S^n\ar@/^20pt/[ruu]^0_{\;}="d"
\ar[r]_<(.3){\S^{n-1}\varepsilon}&S^n\vee S^n\ar[rr]_{\S^{n-1}(k\vee
k)}\ar[uu]_{(1,1)}^{\;}="c"&& S^n_X\vee S^n_X
\ar[uu]_{(1,1)}\ar[rru]^{ev\vee ev}\ar[rr]_<(.35){\S^n\L^n(f
,g)}&&
S^n_Y\ar[r]_{ev}&Y\ar@{=>}"a";"b"_<(.2){(H,0^\vc)}\ar@{=>}"c";"d"^N\ar@{=>}"e";"f"^K}
\end{equation*}
Again by Remark \ref{*} we can remove the nil-track $N$ from (d) and the pasting will still be the same
\begin{equation*}\tag{e}
\xymatrix{&S^n\ar@/^30pt/[rrrr]^0_<(.7){\;}="f"\ar[rr]^{\S^{n-1}k}&&S^n_X\ar[rr]_{ev}^{\;}="e"
&&X\ar@/^20pt/[rdd]^g_{\;}="b"&\\&&&&&X\vee
X\ar[rd]_{(f,g)}^<(.1){\;}="a"\ar[u]^{(1,1)}&
\\S^n
\ar[r]_<(.3){\S^{n-1}\varepsilon}&S^n\vee S^n\ar[rr]_{\S^{n-1}(k\vee
k)}\ar[uu]_{(1,1)}^{\;}="c"&& S^n_X\vee S^n_X
\ar[uu]_{(1,1)}\ar[rru]^{ev\vee ev}\ar[rr]_<(.35){\S^n\L^n(f
,g)}&& S^n_Y\ar[r]_{ev}&Y\ar@{=>}"a";"b"_<(.2){(H,0^\vc)}
\ar@{=>}"e";"f"^K}
\end{equation*}
Pasting some maps and tracks in (e) we obtain
\begin{equation*}\tag{f}
\xymatrix{&&&&&X\ar@/^20pt/[rdd]^g_{\;}="b"&\\&&&&&X\vee
X\ar[rd]_{(f,g)}^<(.1){\;}="a"\ar[u]^{(1,1)}&
\\S^n
\ar[r]_<(.3){\S^{n-1}\varepsilon}&S^n\vee
S^n\ar@/^30pt/[rrrru]^0_{\;}="d"\ar[rr]_{\S^{n-1}(k\vee k)} &&
S^n_X\vee S^n_X \ar[rru]^{ev\vee
ev}^<(.1){\;}="c"\ar[rr]_<(.35){\S^n\L^n(f,g)}&&
S^n_Y\ar[r]_{ev}&Y\ar@{=>}"a";"b"_<(.2){(H,0^\vc)}\ar@{=>}"c";"d"^{K\vee
K} }
\end{equation*}
Again by Remark \ref{*} we can remove $(H,0^\vc)$ from (f) and the pasting of (f) and (g) below coincide
\begin{equation*}\tag{g}
\xymatrix{&&&&&X\vee
X\ar[rd]_{(f,g)}^<(.1){\;}="a"&
\\S^n
\ar[r]_<(.3){\S^{n-1}\varepsilon}&S^n\vee
S^n\ar@/^30pt/[rrrru]^0_{\;}="d"\ar[rr]_{\S^{n-1}(k\vee k)} &&
S^n_X\vee S^n_X \ar[rru]^{ev\vee
ev}^<(.1){\;}="c"\ar[rr]_<(.35){\S^n\L^n(f,g)}&&
S^n_Y\ar[r]_{ev}&Y\ar@{=>}"c";"d"^{K\vee K} }
\end{equation*}
Diagram (g) represents $-f_*[k,K]+g_*[k,K]$, hence equation (3) holds.

The vertical composition of $2$-morphisms $f\st{H}\rr g\st{K}\rr h$ is
preserved by Lemma \ref{general} (3).
The proof of the fact that $\pi_{n,*}$ preserves horizontal
composition is straightforward and it is left to the reader.
\end{proof}

\begin{prop}\label{ki}
The inclusion $\C{cross}_f(n)\subset\C{cross}(n)$ of the full
subcategory of $0$-free objects induces an equivalence of categories
$(n\geq 0)$
$$\C{cross}_f(n)/\!\simeq\;\st{\sim}\To\ho\C{cross}(n),$$
where the homotopy category $\ho$ is obtained by inverting weak equivalences.
\end{prop}

\begin{proof}
For $n=0$ this result is well-known. For $n=1$ this is a consequence
of the fact that $\C{cross}$ has a model category structure where
$0$-free objects are the cofibrant objects, see \cite{htcatg}, and
the homotopy relation derived
from the cylinders on cofibrant objects is given by the tracks defined above. 
In a similar way one obtains the result for $n\geq 2$.
\end{proof}

\section{$k$-Invariants}

Let $K(G,n)$ be the Eilenberg-MacLane space with $\pi_nK(G,n)=G$.
Following Eilenberg-MacLane's notation we write $H^m(G,n,A)$ for
the $m$-dimensional cohomology of the space $K(G,n)$ 
with
coefficients in the abelian group $A$. Here we allow $A$ to be a
$G$-module in case $n=1$. In this case $H^m(G,1,A)=H^m(G,A)$ is the
ordinary cohomology (with local coefficients) of the group $G$.

For any connected $CW$-complex $X$  we write
$$k_n(X)\in H^{n+2}(\pi_nX,n,\pi_{n+1}X)$$
for the first $k$-invariant of the $(n-1)$-connected cover $X\grupo{n}$. Recall that $X\grupo{n}$
is the homotopy fiber of the canonical map from $X$ to its $(n-1)$-type, $X\r P_{n-1}X$, where $P_{n-1}X$ is
a Postnikov section of $X$.

If $n=1$ then $k_1(X)$ is the usual first $k$-invariant of a connected 
$CW$-complex $X$, represented by the crossed module
$$\partial\colon\pi_2(X,X^1)\To\pi_1X^1,$$
determined by the skeletal filtration of $X$, see \cite{3type}.
Otherwise, if $n\geq 2$
$$H^{n+2}(\pi_nX,n,\pi_{n+1}X)=\hom(\Gamma_n\pi_nX,\pi_{n+1}X)$$
and $k_n(X)\colon\Gamma_n\pi_nX\r\pi_{n+1}X$ is induced by the
function $\eta^*\colon\pi_nX\r\pi_{n+1}X$ which sends the homotopy
class of $\alpha\colon S^n\r X$ to the homotopy class of
$\alpha(\S^{n-2}\eta)\colon S^{n+1}\r X$, where $\eta\colon S^3\r
S^2$ is the Hopf map. Compare notation in (\ref{exagam}).

The first secondary homotopy group $\pi_{1,*}X$ is a crossed module,
see Proposition \ref{crosses}. By Proposition \ref{exacta} and
\cite{3type} this crossed module represents an element
$$k(\pi_{1,*}X)\in H^3(\pi_1X,\pi_2X).$$
In general, any crossed module $\partial$ defines a cohomology class
$k(\partial)\in H^3(h_0\partial,h_1\partial),$
see \cite{3type}.

For $n\geq 2$ the $n$-dimensional secondary homotopy group of $X$ defines a homomorphism
$$k(\pi_{n,*}X)\colon\Gamma_n\pi_nX\To\pi_{n+1}X,$$
as follows. Let $k(\pi_{n,*}X)$ be the unique homomorphism fitting into the following commutative diagram
\begin{equation}\label{defica}
\xymatrix{\Gamma_n(\pi_{n,0}X)_\abb\ar@{->>}[d]_{\Gamma_n q }\ar@{^{(}->}[r]&\otimes^2_n(\pi_{n,0}X)_\abb\ar[dd]^\omega\\
\Gamma_n\pi_nX\ar[d]_{k(\pi_{n,*}X)}&\\\pi_{n+1}X\ar@{^{(}->}[r]^\iota&\pi_{n,1}X}
\end{equation}
Here the upper horizontal arrow is the injection in (\ref{exagam}),
and $\iota$ and $q$ appear in Proposition \ref{exacta}.

In general any $0$-free reduced quadratic module $(\omega,\partial)$
defines a homomorphism
$$k(\omega,\partial)\colon\Gamma h_0(\partial,\omega)\To h_1(\omega,\partial),$$
as in (\ref{defica}) and any $0$-free stable quadratic module $(\omega,\partial)$
defines accordingly a homomorphism
$$k(\omega,\partial)\colon h_0(\omega,\partial)\otimes\Z/2\To h_1(\omega,\partial).$$

\begin{thm}\label{sonk}
For any connected $CW$-complex $X$ and any $n\geq 1$ the equality $k_n(X)=k(\pi_{n,*}X)$ holds.
\end{thm}

\begin{proof}
We can suppose without loss of generality that the $1$-skeleton $X^1=\vee_ES^1$ is just a one-point union of $1$-spheres.
One can easily check that $\pi_{1,*}(X,E)$ in Proposition \ref{menos} coincides with
$\partial\colon\pi_2(X,X^1)\r\pi_1X^1$, hence the theorem follows for $n=1$.

We now prove the theorem for $n\geq 2$. Suppose that we have $x\in\pi_{n,0}X$ and we choose $\tilde{x}\colon S^1\r\vee_{\L^nX} S^1$ with $(\pi_1\tilde{x})_\ni(1)=x$. Then
$\omega(\set{x}\otimes\set{x})\in\pi_{n,1}X$ is represented by
\begin{equation*}
\xymatrix{S^n\ar[rr]_{\S^{n-1}\beta}^{\;}="a"\ar@/^30pt/[rr]^0_{\;}="b"&&S^n\vee
S^n\ar@{=>}"a";"b"_{B}\ar[rr]_{\S^{n-1}(\tilde{x},\tilde{x})}&&S^n_X\ar[r]_{ev}&X}
\end{equation*}
The pasting of this is the same as the pasting of
\begin{equation*}
\xymatrix{S^n\ar[rr]_{\S^{n-1}\beta}^{\;}="a"\ar@/^30pt/[rr]^0_{\;}="b"&&S^n\vee
S^n\ar@{=>}"a";"b"_{B}\ar[r]_{(1,1)}&S^n\ar[r]_{\S^{n-1}\tilde{x}}&S^n_X\ar[r]_{ev}&X}
\end{equation*}
By Theorem \ref{propi}
$$\hopf((1,1)B)=-1\in\otimes^2_n\Z=\left\{\begin{array}{cc}\Z,&\text{if $n=2$};\\&\\\Z/2,&\text{if $n\geq 3$}.\end{array}\right.$$
Moreover, $(\pi_1(1,1)\beta)_\ni=0$, therefore by using the definition of $\iota$ in Proposition \ref{exacta}, Theorem \ref{propi}, Remark \ref{conclasi}
and the characterization of $\eta\colon S^3\r S^2$ up to homotopy
as the unique map with Hopf invariant $1$ we get that $\omega(\set{x}\otimes\set{x})=\iota( q (x)(\S^{n-2}\eta))$, hence we are done.
\end{proof}

Let $\C{types}_n^1$ be the category of pointed $(n-1)$-connected $CW$-complexes $X$ with $\pi_m(X,x_0)=0$ for all $m\geq n+2$ and all $x_0\in X$.

\begin{prop}\label{loctite}
The functor $\pi_{n,*}\colon\C{types}_n^1\r\C{cross}(n)$ induces an equivalence of categories 
$(n\geq 0)$
$$\pi_{n,*}\colon\ho\C{types}_n^1\st{\sim}\To\ho\C{cross}(n),$$
where the homotopy category $\ho$ is obtained by localizing with respect to weak equivalences.
\end{prop}

For the proof of Proposition \ref{loctite} we recall the following functors.

Let $\C{CW}_n$ be the category of $CW$-complexes $X$ with trivial
$(n-1)$-skeleton $X^{n-1}=*$ and cellular maps. There is a
``cellular'' functor
\begin{equation}\label{cell}
P_{n+1}\sigma\colon\C{CW}_n/\!\simeq\,\To\C{cross}_f(n)/\!\simeq.
\end{equation}
If $n=1$ this functor sends a $CW$-complex $X$ to the crossed module
$$\partial\colon\pi_2(X,X^1)\r\pi_1X^1$$
given by the boundary operator in the long exact sequence of
homotopy groups, see \cite{ctagIII} and \cite{3type}. If
$n\geq 2$ the the reduced (stable if $n\geq 3$) quadratic module
$P_n\sigma(X)$ is the truncation of the totally free quadratic
complex $\sigma(X)$ defined in \cite{ch4c} IV.C,
$$\otimes^2C_n(X)\st{\omega}\To\sigma_{n+1}(X)/d(\sigma_{n+2}(X))\st{\partial}\To\sigma_n(X),$$
compare \cite{ch4c} IV.10.4 
and \cite{cosqc} 4.

\begin{prop}\label{concel}
The functor $P_{n+1}\sigma$ in (\ref{cell}) is naturally isomorphic
to $$\pi_{n,*}\colon\C{CW}_n/\!\simeq\,\r\C{cross}_f(n)/\!\simeq$$ for
all $n\geq 1$.
\end{prop}

\begin{proof}
If $X$ is $(n-1)$-reduced then $X^n=\vee_ES^n$ for some pointed set $E$. The inclusion of spheres in the wedge
$X^n\subset X$
determines a pointed inclusion $E\subset\L^nX$. 
One can easily check that $P_{n+1}\sigma(X)$ is isomorphic to $\pi_{n,*}(X,E)$ in Proposition \ref{menos}. 
Now the natural isomorphism in the statement is given by the weak equivalence $P_{n+1}\sigma(X)\cong\pi_{n,*}(X,E)
\st{\sim}\r\pi_{n,*}X$ in Proposition \ref{menos}. Compare Proposition \ref{ki}.

\end{proof}

\begin{proof}[Proof of Proposition \ref{loctite}]
For $n=0$ this is a well-known result. For $n\geq 1$ this follows from Propositions \ref{ki}, \ref{concel} and the fact that
$P_{n+1}\sigma$ in (\ref{cell}) does induce an equivalence of categories $P_{n+1}\sigma\colon\ho\C{types}_n^1\r\ho\C{cross}(n)$.
This is shown in \cite{ch4c} III.8.2 for $n=1$. For $n=2$ the proof follows as in the case $n=1$, this case is considered  even
in the non-simply connected case in \cite{ch4c} IV.10.1. The case $n\geq 3$ can be easily proved along the lines
of the $n=1$ and $n=2$ cases,
i.e. by using \cite{ch4c} III.8.5, III.8.8 and IV.C.14. 
\end{proof}


\bibliographystyle{amsalpha}
\bibliography{Fernando}
\end{document}